\documentclass{siamart0516}
\title{Stability estimates for linearized near-field phase retrieval in X-ray phase contrast imaging}

\usepackage[utf8]{inputenc}
\usepackage{cite}
\usepackage{amstext,amsmath,amsfonts,amssymb,amsopn,ifthen,fancyhdr,lastpage,mathrsfs,color,listings,dsfont}
\usepackage{graphicx}
\usepackage{xcolor}
\usepackage{subfig}

\usepackage[textwidth = 2cm,textsize=small]{todonotes}

\newcommand{\Fref}[1]{Fig.~\ref{#1}}
\newcommand{\eref}[1]{(\ref{#1})}

\newcommand{\sC}{\mathscr{C}}

\newcommand{\cD}{\mathcal{D}}

\newcommand{\cF}{\mathcal{F}}

\newcommand{\cN}{\mathcal{N}}
\newcommand{\cO}{\mathcal{O}}

\newcommand{\mN}{\mathbb{N}}

\newcommand{\mR}{\mathbb{R}}

\newcommand{\mS}{\mathbb{S}}

 \newcommand{\bj}{\ensuremath{\boldsymbol{j}}}
 \newcommand{\bx}{\ensuremath{\boldsymbol{x}}}
 
  \newcommand{\bhhat}{\ensuremath{\boldsymbol{\hat h}}}
 \newcommand{\bk}{\ensuremath{\boldsymbol{k}}}

\newcommand{\bxi}{\ensuremath{\boldsymbol{\xi}}}
 
 \newcommand{\bn}{\ensuremath{\boldsymbol{n}}}

\newcommand{\bepsilon}{\ensuremath{\boldsymbol{\epsilon}}}

\newcommand{\btheta}{\ensuremath{\boldsymbol{\theta}}}

\newcommand{\Text}[1]{\text{\textnormal{#1}}}

\newcommand{\I}{\Text{i}}
\newcommand{\E}{\Text{e}}

\newcommand{\MTEXT}[1]{\;\;\;\;\;\text{#1}\;\;\;\;\;}
\newcommand{\ip}[2]{\left\langle #1, #2\right\rangle} % Skalarprodukt
\newcommand{\D}{\Text{d}}
\newcommand{\norm}[1]{\left\| #1 \right\|} % Norm
\newcommand{\bnorm}[1]{\big\| #1 \big\|} % Norm
\newcommand{\closure}[2][3]{%
  {}\mkern#1mu\overline{\mkern-#1mu#2}}

\newcommand{\tNF}{\mathfrak{f}}%{\tilde N_{\textup{F}}}
\newcommand{\tNFa}[1]{\tNF_{#1}}%{\tilde N_{\textup{F}, #1}}
\newcommand{\Fres}{\lower0.9ex\hbox{$\mathchar'26$}\mkern-8mu \tNF}

\newcommand{\nF}{n_{\tNF}}%{n_{\textup{F}}}
\newcommand{\mF}{m_{\tNF}}%{m_{\textup{F}}}
\newcommand{\OmegaF}{\Omega_{\tNF}}%{\Omega_{\textup{F}}}
\newcommand{\IF}{I_{\tNF}}%{I_{\textup{F}}}
\newcommand{\FF}{\mathcal{F}_{\tNF}}%{\mathcal{F}_{\textup{F}}}

\newcommand{\cc}[1]{\overline{#1}}

\newcommand{\supp}{\textup{supp}}

\newcommand{\comp}[1]{#1^\textup{c}}

\newcommand{\secref}[1]{$\S$\ref{#1}}

\newtheorem{MyIP}{Inverse Problem}
\crefname{MyIP}{Inverse Problem}{Inverse Problems}

\newcommand{\TheTitle}{Stability estimates for linearized near-field phase retrieval in X-ray phase contrast imaging} 
\newcommand{\TheAuthors}{S. Maretzke, T. Hohage}

\headers{Stability of X-ray phase contrast imaging}{\TheAuthors}

\title{{\TheTitle}\thanks{\textbf{Published as:} SIAM J.\ Appl.\ Math. (2017), 77(2), 384--408. DOI:10.1137/16M1086170}
\funding{This work was funded Deutsche Forschungsgemeinschaft DFG through Project C02 
of SFB 755 - Nanoscale
Photonic Imaging}}

\author{
  Simon Maretzke\thanks{Institute for Numerical and Applied Mathematics, University of G\"ottingen, Lotzestrasse 16-18, 37083 G\"ottingen, Germany
    (\email{s.maretzke@math.uni-goettingen.de}).} %, \url{http://www.imag.com/\string~ddoe/}
  \and
  Thorsten Hohage\thanks{Institute for Numerical and Applied Mathematics, University of G\"ottingen, Lotzestrasse 16-18, 37083 G\"ottingen, Germany (\email{hohage@math.uni-goettingen.de}).}
}

\begin{document}
\maketitle

\begin{abstract}
Propagation-based X-ray phase contrast enables nanoscale imaging of 
biological tissue by probing not only the attenuation, but also 
the real part of the refractive index of the sample. Since only intensities 
of diffracted waves can be measured, the main mathematical challenge consists 
in a phase-retrieval problem in the near-field regime. We treat an
often used linearized version of this problem known as 
contract transfer function model.
Surprisingly, this inverse problem turns out to be well-posed assuming 
only a compact support of the imaged object. Moreover, we establish bounds on the 
Lipschitz stability constant.
In general this constant grows exponentially with the Fresnel number of the imaging setup. 
However, both for homogeneous objects, characterized by  a fixed ratio 
of the induced refractive phase shifts and attenuation,
and in the case of measurements at two distances, a much more favorable algebraic
dependence on the Fresnel number can be shown. In some cases we 
establish order optimality of our estimates. 
\end{abstract}

% Uncomment for PACS numbers
%\pacs{02.30.Nw, 02.30.Zz, 42.30.Rx, 42.30.Wb}



\section{Introduction} \label{S1}

Over the past two decades, the dramatic increase in coherence and brightness 
of large-scale X-ray sources, such as third generation synchrotrons and 
free-electron-lasers, has paved the way for X-ray phase contrast imaging 
\cite{Nugent2010coherent}. Classical X-ray radiography is limited to measuring 
the attenuation experienced by radiation traversing the probed object. 
Writing the refractive index in the X-ray physics notation 
$n = 1-\delta+\I \beta$ with $0\leq \beta,\delta\ll 1$, this amounts to imaging $\beta$. 
Phase contrast techniques additionally probe the real-valued decrement
$\delta$ of $n$, which induces
phase shifts in the transmitted X-ray wave field. 
This enables imaging of biological cells and other micro-scale 
light-element specimen, for which $\beta\ll \delta$ holds in the hard X-ray regime
\cite{Wilkins1996,Paganin1998,Cloetens1999,Mayo2002quantitative}. Owing to the small wavelength of X-rays, nano-scale spatial resolutions can be achieved as has been demonstrated down to 20 nanometers \cite{Bartels2015}. Moreover, phase contrast imaging can be combined with tomography, capable of resolving the refractive index of an unknown object in 3D \cite{Jonas2004TwoMeasUniquePhaseRetr,Bartels2012,Kostenko2013_AllAtOncePCTWithTV,Krenkel2014BCAandCTF,Ruhlandt2016_RadonPCICommute,MaretzkeEtAl2016OptExpr}.

Unfortunately, the refractive phase shifts of the X-ray field 
% in the exit plane $z=0$ (see Fig.~\ref{figure1}) 
cannot be observed directly by common CCD detectors due to their physical limitation to measuring wave intensities, i.e.\ the squared modulus of the wave field. 
In propagation-based phase contrast imaging, also known as \emph{inline holography}, the required phase-sensitivity is achieved simply by free-space propagation without any optical elements: if the detector is placed in some finite distance down-stream of the sample, the imprinted phase shifts in the object's exit plane ($z=0$ in \cref{figure1}(b))  are partially encoded into measurable intensities by \emph{diffraction}, i.e.\ self-interference of the wave field.
% To recover the field in the exit plane we have to solve a \emph{phase retrieval} problem. 
We assume that the diffraction pattern or \emph{hologram} is recorded in the optical near-field of the sample so that propagation is described by the Fresnel propagator \cite{PaganinXRay}. In particular, we do not consider the corresponding far-field setup (coherent diffactive imaging, see  e.g.\ \cite{Quiney2010CDIReview,Miao2015XrayCDIReview}), where the data is given by Fourier magnitudes.

In this work, we are thus concerned with the reconstruction of the (generally complex-valued) wave-field perturbation $h$ induced by the object from measured near-field intensities $I$. As this implicitly amounts to recovering the lost phase information in the data, i.e.\ to solving a \emph{phase retrieval} problem, the question immediately arises whether the image recovery is actually \emph{unique} and \emph{stable}. Indeed, it is commonly argued \cite{Jonas2004TwoMeasUniquePhaseRetr,Nugent2007TwoPlanesPhaseVortex,Burvall2011TwoPlanes} that diffraction patterns from at least two different sample-detector distances are required for a unique reconstruction of the imprinted phase shifts and attenuation. Assuming a \emph{support constraint}, however, i.e.\ under the often physically reasonable assumption that the image $h$ is non-zero only in compact subdomain of the field of view, 
we could show uniqueness of the reconstruction from a \emph{single} hologram in \cite{Maretzke2015IP} - even if $h$ is complex-valued (see also 
\cite{klibanov:14} for a uniqueness result for the Helmholtz equation 
with real-valued $n$ from phaseless near-field data in an interval of frequencies). 
In \cite{MaretzkeEtAl2016OptExpr}, such a reconstruction of a compactly supported complex image is demonstrated for simulated and experimental data, which turns out to be feasible, yet susceptible to low-frequency artifacts. Support constraints have also been found to stabilize image reconstruction in simpler settings where the probed object can be assumed to be completely non-absorbing \cite{Giewekemeyer2011CellSketch,Bartels2012}.
% Phase-sensitivity can be achieved for instance by interferometry\cite{Bonse1965interferometric,Momose1995interferometric,Wilkins1996interferometric} but even without any optical elements just by free-space propagation: 

These observations call for a better understanding of the \emph{stability} of the considered near-field phase retrieval problem, which is the goal of this paper. 
We do so within a linearization of the relation between image $h$ and the resulting intensity data $I$, valid for sufficiently ``small'' $h$, i.e.\ for weakly interacting objects 
similarly as in a recent stability analysis of domain reconstructions in phaseless inverse scattering \cite{AmmariEtAl2016_StabAnalysisLinInvScat}.
The linearization is known as \emph{contrast transfer function} model \cite{Guigay1977CTF,Turner2004FormulaWeakAbsSlowlyVarPhase} and frequently applied in X-ray phase contrast imaging \cite{Cloetens1999,Gureyev2004CTF,HoffmannEtAl2011_CTFCriticality,LangerEtAl2012_PCT_CTF,Krenkel2014BCAandCTF}. In this work, we analyze the arising linear forward operator $T$
%: h \mapsto  Th \approx I -I_0$ 
under the assumption that $h$ has compact support. 
%support- and additional constraints on $h$ for the settings of one or two recorded holograms. 
We prove that %, under fairly mild constraints,
the associated inverse problem is not only unique but even \emph{well-posed} in the sense that the recovered image $h$ depends continuously on the measured data $Th$, i.e.\ finite data errors lead to 
bounded deviations in the recovered $h$. This result, which is quite surprising 
for an inverse problem with remote measurements, is achieved by relating the setting
to a reconstruction from incomplete Fourier data. 
By the same technique we also derive explicit \emph{stability estimates}, bounding the reconstruction error that results from a given data noise level in terms of the dimensionless Fresnel number.
%, which is inverse proportional to the distance between the sample and the detector. 

In general, we find that the stability constant decays exponentially with the Fresnel number $\tNF$ and thus hardly gives any useful stability bounds for many experimental 
X-ray phase contrast setups. However, we establish much more favorable $\cO(\tNF^{-1})$- and even $\cO(\tNF^{-1/2})$-decay rates  % of the stability constant with the inverse Fresnel number 
in two relevant situations: 
The first assumes proportionality of the real and the imaginary part of the 
image $h$, which occurs e.g.\ for single-material and 
non-absorbing samples. The second situation concerns general objects, but 
two measurements at different distances. It is well-known that the forward operators 
in these cases are Fourier multipliers with so-called 
contrast transfer functions (CTF), and the zeros of the CTFs are responsible 
for general ill-posedness of the inverse problem. % without support constraints.
Our analysis exploits the regularizing effect on these zeros of the smoothness in  
Fourier space that results from the assumed support constraints. 
%As the forward operator $S_{\alpha}$ in \eqref{eq:fwModel2} arises from $T$ by a mere \emph{restriction} to homogeneous images, well-posedness and stability of \cref{ip1} directly carry over to \cref{ip2}. However, it should be emphasized that the proven \emph{exponential} rates for the constant $C_{\textup{IP1}}(\Omega,\tNF) \gtrsim \exp(-\tNF / 8)$ hardly guarantee stability in a practical sense for Fresnel numbers $\Fres \geq 100$, as encountered in typical near-field imaging experiments. We have seen in \secref{SS:StabResults} that \cref{ip2} is in general much less ill-posed than \cref{ip1}. In this spirit, \cref{thm:IP2StabRes} establishes improved \emph{algebraic} rates for the stability constant $C_{\textup{IP2}}(\Omega,\tNF) \gtrsim \tNF^{-1}$ with the Fresnel number $\tNF$. The aim of the following analysis is to prove this stability result.

The remainder of this paper is organized as follows: in \secref{S2}, the mathematical model of X-ray phase contrast imaging is introduced and the considered inverse reconstruction problems, corresponding to different constraints and measurement setups, are motivated. Our principal stability results are stated in \secref{SS:StabResults}. \secref{S3}, \secref{S4} and \secref{S5} contain the analysis for each of the inverse problems, including the proofs of the main results. In \secref{S6}, we discuss implications of our findings and possible extensions.

\vspace{1em}

\section{Imaging problems and main results} \label{S2}

\subsection{Physical model} \label{SS2.1}

\begin{figure}[hbt!]
 	 \centering
 	 \includegraphics[width=.8\textwidth]{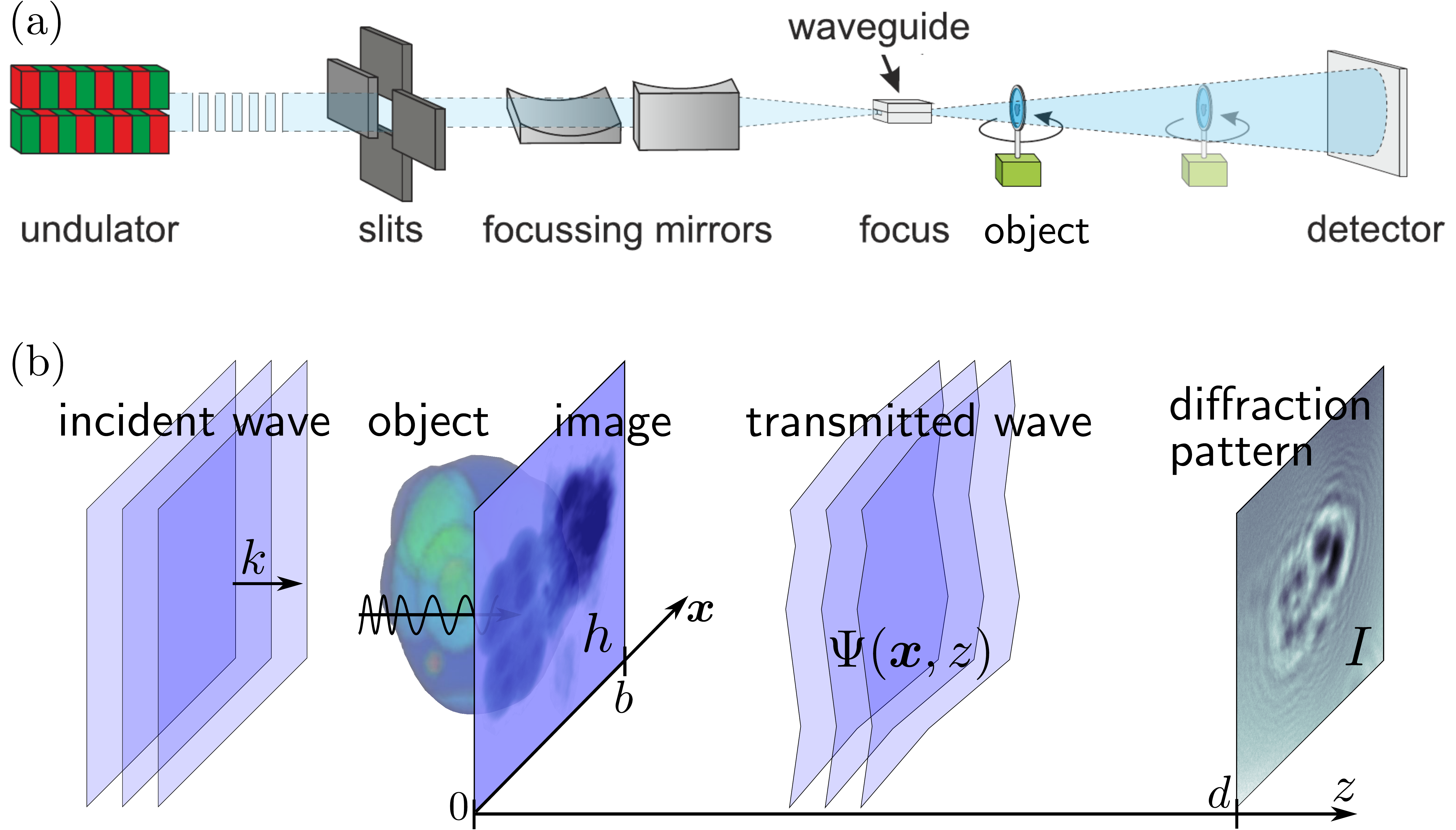}
 	 \caption{Setup of propagation-based X-ray phase contrast imaging (inline holography).\\
 	 (a) Sketch of an experimental realization \cite{MaretzkeEtAl2016OptExpr} (GINIX \cite{Salditt2015GINIX} at P10-beamline, DESY) \\
 	 (b) Physical model \cite{Maretzke2015IP}: incident plane waves are scattered by an unknown sample, imprinting a phase- and absorption image $h =-\I \phi - \mu$ upon the transmitted wave fronts $\Psi(\cdot,0)$. The resulting near-field diffraction pattern (hologram) $I= |\Psi(\cdot, d)|^2 $ is recorded at some distance $z = d$ behind the object. (plotted experimental data: hologram $I$ and reconstructed image $h$ of \emph{d.\ radiodurans} bacteria). \label{figure1}}
 	\end{figure}

An exemplary experimental setup for X-ray phase contrast imaging on a third generation synchrotron source (GINIX setup \cite{Salditt2015GINIX} at P10-beamline, DESY) is shown in \cref{figure1}(a).
We describe this imaging system by a standard wave-optical model as schematically visualized in  \cref{figure1}(b)  \cite{Pogany1997noninterferometric,Cloetens1999Diss}: an unknown sample is illuminated by an incident plane electromagnetic wave $\Psi_{\textup i}(\bx, z) = \exp(\I k z)$, where $\bx \in \mR^2$ and $z\in \mR$ denote the lateral- and axial coordinates, respectively. By scattering interaction, object information is encoded as a perturbation of the wave field $\Psi = \Psi_{\textup i} + \Psi_{\textup s}$ within the exit plane $z =0$ of the probed sample. A detector measures the resulting near-field diffraction pattern (or hologram), given by the intensity $I(\bx) = |\Psi(\bx,d)|^2$ of the propagated wave fronts in some plane at finite distance $d>0$ behind the object. The phase of the complex-valued field $\Psi$ cannot be observed directly, yet \emph{diffraction} partially encodes phase variations in the exit plane $z= 0$ into measurable intensities $I$ at the detector.
 	
In general, the above physical model is governed by the Helmholtz equation $\Delta \Psi + n^2 k^2 \Psi = 0$ in $\mR^3$, where the object gives rise to a spatially varying refractive index $n = 1 - \delta + \I \beta$.
In the hard X-ray regime of very large wavenumbers $k$, typical samples such as  biological cells are often sufficiently thin and weakly interacting for the scattering  to be well-approximated by geometrical optics \cite{PaganinXRay,Jonas2004TwoMeasUniquePhaseRetr}. Within this approximation, the wave fronts $\Psi_z:= \Psi(\cdot, z)$ in the exit plane are given by
\begin{equation}
  \Psi_{0} =   \exp(h)  \MTEXT{with} h = -\I \phi - \mu =- \I k \int_{\mR} \big( \delta - \I \beta \big)  \D z  \label{eq:2.1-1}
\end{equation}
Accordingly, the perturbed wave yields line integrals over $\delta$ and $ \beta$ along the incident $z$-direction, corresponding to a \emph{projection image} of the sample in the form of phase shifts $\phi$ and attenuation $\mu$. In particular, \eqref{eq:2.1-1} implies that the imprinted image $h$ satisfies a \emph{support constraint} $\supp(h) \subset \Omega$  whenever the object is laterally finite, i.e.\ if $\delta(\bx,z) = \beta(\bx,z) = 0$ for all $z \in \mR, \bx \in \Omega  $ outside some bounded domain $\Omega \subset \mR^2$.

If $|\bxi|\ll k$ for all relevant spatial frequencies $\bxi$ of $h$, 
the total wave field will be of the form $\Psi(\bx,z)=\E^{\I kz}\tilde{\Psi}(\bx,z)$ with a slowly 
varying envelope $\tilde{\Psi}$ such that 
$\E^{-\I kz} (\partial_z^2+k^2)\Psi % = (-\I\partial_z+k)(\I\partial_z+k)\Psi
\approx  2\I k \partial_z \tilde \Psi$ by neglecting $\partial_z^2\tilde\Psi$.  
This yields the one-way-, Schr\"odinger- or paraxial approximation
\[%begin{equation}\label{eq:Schroedinger}
2k\I\partial_z\tilde{\Psi}+\Delta_{\bx}\tilde{\Psi}\approx 0\qquad 
%\mbox{for}\qquad \tilde{\Psi}(\bx,z):=\E^{-\I\kappa z}\Psi(\bx,z)
\]%end{equation}
to the Helmholtz equation.  
Within this commonly used model, the free-space propagation 
of the wave fronts to the detector is described by the \emph{Fresnel propagator} 
\cite{PaganinXRay}:
\begin{eqnarray}
 \cD ( \Psi_0 ) &:=& \exp(-\I kd)\Psi_d =   \cF^{-1}\left( \mF \cdot \cF(\Psi_0) \right) 
\qquad 
 \mF(\bxi )  := \exp \left(  \frac{  -\I |\bxi|^2 }{2\tNF}  \right)  \label{eq:FresnelProp}
\end{eqnarray}
Here, $\cF$ is the Fourier transform and 
\begin{equation}\label{eq:defi_Fresnel}
\tNF := \frac{kb^2}{d}\qquad \mbox{or}\qquad 
\Fres := \frac{\tNF}{2\pi} 
\end{equation}
denote the dimensionless \emph{Fresnel number} of the setup. 
$b$ is a physical length that corresponds to length 1 in dimensionless coordinates 
and will be chosen as the support diameter of the image $h$, see \cref{figure1}(b). 
Usually $\Fres$ is referred to as Fresnel number, but this convention would lead to an 
abundance of $2\pi$ factors in our computations. Therefore, we will mostly use 
$\tNF$ with the notation \eqref{eq:defi_Fresnel} chosen in analogy to 
Planck's constant. $\tNF$ governs the impact of diffraction in the imaging setup, where smaller values correspond to stronger diffractive distortion of the propagating wave field. 
Typical values in experimental X-ray phase contrast setups are in the range $10 \leq \Fres \leq 1000$.
%\[
%\tNF := \frac{kb^2}{d}
%\] 
%denotes the dimensionless (modified) \emph{Fresnel number} of the setup. $b$ is a physical length that corresponds to length 1 in dimensionless coordinates. Usually 
%$\Fres$ is referred to as Fresnel number, but this convention would lead to an 
%abundance of $2\pi$ factors in our computations.
%% The scaling does not 
%% matter for asymptotic analysis, but if numerical values are concerned we will 
%% always specify the values of $\Fres$. 
%$\tNF$ governs the impact of diffraction in the imaging setup, where smaller values correspond to stronger diffractive distortion of the propagating wave field.  Note that we take $b$ as the support diameter of the image $h$, see \cref{figure1}(b). Typical values of the
%standard Fresnel number $\Fres$ in experimental X-ray phase contrast setups are then 
%in the range $10 \leq \Fres \leq 1000$.

By combination of \eref{eq:2.1-1} and \eref{eq:FresnelProp}, we find that the unknown object image $h$ is related to the observable intensity data $I = |\Psi_d|^2$ by the nonlinear forward operator
\begin{equation}
 I = F(h) := \big|  \cD\big( \exp(h)  \big)  \big|^2 \MTEXT{with} h =  -\I\phi- \mu \label{eq:fwModelNL}
\end{equation}
We note that similar models apply to imaging with \emph{electrons} \cite{Erickson1971_CTF_TEM,Wade1992_CTF_TEM,Lichte2007_ElectronHolography} owing to the mathematical equivalence of the time-independent Schr\"odinger equation and the Helmholtz equation.

\vspace{.5em}

\subsection{Weak object limit and principal inverse problem} \label{SS2.2}

By \eqref{eq:fwModelNL}, the image $h$ is in general complex-valued, whereas the intensity data $I$ is real-valued. This suggests that the data is insufficient for unique and stable recovery of $h$ \cite{Nugent2007TwoPlanesPhaseVortex,Burvall2011TwoPlanes}. We analyze this question of ill-posedness within the   commonly used \emph{weak-object-approximation} \cite{Pogany1997noninterferometric,Cloetens1999,Paganin2002simultaneous,Gureyev2004CTF}: in the case of weak absorption $\mu \ll 1$ and slowly varying phase shifts $\phi$, nonlinear terms in $h$ can be neglected in \eref{eq:fwModelNL} \cite{Turner2004FormulaWeakAbsSlowlyVarPhase}, giving
\begin{eqnarray}
 F(h) = 1 + Th + \cO(h^2) \MTEXT{with} Th &:= 2\Re\big(\cD(h) \big).  \label{eq:fwModel1}% = \cD(h) + \closure{\cD(h)}  
 %&:= 2\Re\big(\cD(h) \big) = \cD(h) + \closure{\cD(h)}.  \label{eq:2.2-1-b} 
\end{eqnarray}
Here, $\Re$ denotes the pointwise real-part and we have used that $\cD(1) = 1$.
% Rigorous error estimates for the linearization \eqref{eq:fwModel1} are given in \secref{SS5.1}. 
  Note that $\cD$ is unitary on $L^2(\mR^m)$, so $T$ defines a bounded $\mR$-linear operator on $L^2(\mR^m)$.
 Rather than by \eqref{eq:fwModel1}, $T$ is more commonly written in terms of the phase shifts $\phi$ and absorption $\mu$ via sinusoidal \emph{contrast transfer functions} (CTF) \cite{Guigay1977CTF}:
  \begin{equation}
  T (-\I\phi- \mu) =  2 \cF^{-1} \left( \sin \left( \frac{|\bxi|^2 }{2\tNF} \right) \cF( \phi )  - \cos \left( \frac{|\bxi|^2 }{2\tNF} \right) \cF(  \mu) \right).  \label{eq:fwModelCTF}
 \end{equation}

Although the physical model of \secref{SS2.1} leads to two-dimensional images $h$ and holograms $I$, we will study the operator $T$ in a more general $\mR^m$-setting. This might allow application of our results to situations described by (quasi-) 1D-models and to \emph{phase contrast tomography}, which can be interpreted as a \emph{3D}-imaging modality \cite{Kostenko2013_AllAtOncePCTWithTV,Ruhlandt2016_RadonPCICommute}. As physically motivated in \secref{SS2.1}, we impose support constraints by assuming
\begin{equation}
h \in L^2_\Omega := \{ h \in L^2(\mR^m): h|_{\mR^m \setminus \Omega} = 0 \} \MTEXT{for some} \Omega \subset \mR^m.
\end{equation} 
% for some support-domain $\Omega\subset \mR^m$. 
Moreover, we denote by $\norm{h} := (\int_{\mR^m} |h|^2 \; \D\bx )^{1/2}$ the standard $L^2$-norm in $\mR^m$. The principal image reconstruction problem of this work then reads as follows:

\vspace{.5em}
\begin{MyIP}[Phase contrast imaging of weak objects]\label{ip1}
 For a given support $\Omega \subset \mR^m$, recover a complex-valued image  $h \in L^2_\Omega$ from noisy intensity data 
 \begin{equation*} 
  I^{\bepsilon} = 1 + T  h + \bepsilon \MTEXT{with} \norm{\bepsilon}  \leq \epsilon.
 \end{equation*}
 %with deterministic errors $\bepsilon$ s.t. $\norm{\bepsilon}_{ L^2(\mR^m)} \leq \delta$.
\end{MyIP}
\vspace{.5em}
% 
% \noindent Here, $A$ denotes the set of admissible solutions, defined by \emph{a priori} knowledge on the images to be recovered.

\subsection{Homogeneous and non-absorbing objects} \label{SS:HomObj}

It is often legitimate to assume that the object is \emph{homogeneous} in the sense that phase shifts $\phi$ and attenuation $\mu$ are proportional, i.e.\
\begin{equation}
  h =  -  \mu - \I \phi = -\I \E^{-\I \alpha} \varphi \label{eq:SingleMat} 
\end{equation}
%  This assumption is satisfied in particular if the object is composed of a single material. The constraint can be incorporated by setting
  for some $\alpha \in [0;\pi)$ and a  \emph{real-valued} function $\varphi \in L^2(\mR^m)$. Note that this includes the special case $\alpha = 0$ which corresponds to $\mu = 0$ and thus to a purely phase shifting, i.e.\ non-absorbing object, providing an excellent model for hard X-ray imaging of light-element samples. By plugging \eqref{eq:SingleMat} into \eqref{eq:fwModelCTF}  and rearranging by trigonometric identities,  we obtain 
	a forward operator incorporating the \emph{homogeneity constraint}:
\begin{equation}
\begin{array}{rl}
S_{\alpha} : L^2(\mR^m )  &\!\!\!\!\to L^2 (\mR^m),\\
 \varphi  &\!\!\!\!\mapsto  T (-\I \E^{-\I \alpha} \varphi)  = 2 \cF^{-1} \left( s_\alpha \cdot \cF (\varphi) \right) 
\end{array}\qquad 
s_\alpha(\bxi ) := \sin \left( \frac{ |\bxi|^2 }{2\tNF} + \alpha \right)\label{eq:fwModel2}
\end{equation}
Accordingly, the forward model reduces to a multiplication with the contrast transfer function $s_\alpha$ (CTF) in Fourier space \cite{Guigay1977CTF,Turner2004FormulaWeakAbsSlowlyVarPhase}. %, which is often exploited in direct reconstruction methods.
This makes the inversion of $S_{\alpha}$ significantly easier than that of $T$, which is why we state it as a second inverse problem: 
 
\vspace{.5em}
\begin{MyIP}[Phase contrast imaging of weak homogeneous objects]\label{ip2}
 For given $\Omega \subset \mR^m$, recover a real-valued image  $\varphi \in L^2_\Omega$ from noisy intensity data
 \begin{equation*} 
  I^{\bepsilon} = 1 + S_{\alpha}  \varphi + \bepsilon \MTEXT{with} \norm{\bepsilon}  \leq \epsilon.
 \end{equation*}
\end{MyIP}
\vspace{.5em}

\subsection{Stability estimates} \label{SS:StabResults}

The statement of \cref{ip1,ip2} immediately raises the question whether these are uniquely solvable and whether the solution is \emph{stable} with respect to noise $\bepsilon$.
In order to illustrate the significance of this problem, we first recall some well-known facts on the derived inverse problems without assuming a support constraint, i.e.\ for $\Omega = \mR^m$: as the null-space $\cN(T) = \{\I \cD^{-1}(f): f \in L^2(\mR^m)\text{ real-valued}\}$ of the forward map $T: L^2(\mR^m) \to L^2(\mR^m)$ is 
huge, \cref{ip1} is heavily non-unique in this setting. The operator $S_{\alpha}$, on the other hand, is indeed injective so that \cref{ip2} is uniquely solvable. However, the inversion of $S_{\alpha}$ is \emph{ill-posed} since noise in Fourier-frequencies near the zeros of the CTF $s_\alpha$ is amplified by arbitrary factors in the reconstruction.% The same holds true for \cref{ip3}

To our great surprise, imposing a support constraint with bounded $\Omega$ does not 
only rule out non-uniqueness (as proven in \cite{Maretzke2015IP}), but even turns \cref{ip1} into a \emph{well-posed} problem: every admissible image $h \in L^2_\Omega$ gives rise to finite \emph{contrast} $\norm{Th} \geq C_{\textup{IP1}} \norm{h}$ in the observable data with some lower bound $C_{\textup{IP1}} > 0$:

\vspace{.5em}
\begin{theorem}[Well-posedness and stability estimate for \cref{ip1}] \label{thm:IP1StabRes} 
Let the support-domain
 $\Omega$ be given by a stripe of width $1$, without loss of generality \ $\Omega := [-1/2; 1/2] \times \mR^{m-1}$. Then there exists a constant $C_{\textup{IP1}}(\Omega,\tNF) > 0$ such that
 \begin{equation}
  \norm{T  h}  \geq C_{\textup{IP1}}(\Omega,\tNF)  \norm{h} \MTEXT{for all} h \in L^2_\Omega, \label{eq:IP1StabRes-1}
 \end{equation}
 i.e.\ \cref{ip1} is well-posed. The stability constant satisfies the estimate
 \begin{equation}
  C_{\textup{IP1}}(\Omega,\tNF)\geq  (2\pi \tNF)^{\frac 1 4} \bigg( 1 - \frac{3}{8 \tNF } + \cO\Big(\tNF ^{-2} \Big) \bigg)   \exp\left(-\tNF / 8 \right). \label{eq:IP1StabRes-2}
 \end{equation}
 \end{theorem}
\vspace{.5em}

  \cref{thm:IP1StabRes} is proven in \secref{S3} along with a characterization of the least stable modes, i.e.\ of the images $h$ that induce least contrast under $T$. Notably, \eqref{eq:IP1StabRes-1} implies $\norm{h} \leq C_{\textup{IP1}}(\Omega,\tNF)^{-1}  \norm{Th}$, ensuring finite amplification of data errors  $\leq C_{\textup{IP1}}(\Omega,\tNF)^{-1} \norm{\bepsilon}$ upon inversion of $T$ and thus \emph{stability} of the reconstruction of $h$ from $I^{\bepsilon}$. By \eqref{eq:IP1StabRes-2}, however, the constant $C_{\textup{IP1}}(\Omega,\tNF)$ decays (nearly) exponentially with increasing $\tNF$ so that \cref{thm:IP1StabRes}  hardly guarantees stability in any practical sense for Fresnel numbers $\Fres \geq 100$. Fortunately, the stability estimate can be  improved to algebraic decay  with $\tNF$ in the case of \cref{ip2}, as shown in \secref{S4}:

\vspace{.5em}
\begin{theorem}[Well-posedness and stability estimate for \cref{ip2}] \label{thm:IP2StabRes}
Let the support-domain $\Omega := \{ \bx \in \mR ^m : |\bx | \leq \frac 1 2\} $ be a ball of diameter $1$. Then
the stability constant $C_{\textup{IP2}}(\Omega,\tNF, \alpha):= \inf_{\varphi\in L^2_\Omega, \norm{\varphi} = 1} \norm{S_\alpha \varphi}$
of \cref{ip2} is bounded by
 \begin{equation}
  C_{\textup{IP2}}(\Omega,\tNF, \alpha)  \geq \max \left\{ \min \left\{  c_{1}  ,  c_{2} \tNF^{-1} \right\}, \min \left\{  c_{3} \alpha  ,  c_{4} \tNF^{-\frac 1 2} \right\} \right\} \label{eq:IP2StabRes-2}
 \end{equation}
for some constants $c_j > 0$ that depend only on the dimension $m$. 
In particular, $C_{\textup{IP2}}(\Omega,\tNF, \alpha)=\mathcal{O}(\tNF^{-1})$ 
for $\alpha=0$ and  $C_{\textup{IP2}}(\Omega,\tNF, \alpha)=\mathcal{O}(\tNF^{-1/2})$ 
for $\alpha>0$ as $\tNF\to\infty$. 
\end{theorem}
\vspace{.5em}

% In \secref{S5}, we show that stability of \cref{ip3} is closely related to the forward operator $S$ of \cref{ip2} for $\alpha = 0$. Hence, the same improvements of the stability constant compared to \cref{thm:IP1StabRes} hold in this setting:
% 
% 

We recall that the physical lengthscale $b$ underlying to the Fresnel number $\tNF$ in \cref{thm:IP1StabRes,thm:IP2StabRes} is the diameter of the support-domain $\Omega$ as the latter is taken to be unit length. Accordingly, the resulting stability constants $C_{\textup{IP1}}, C_{\textup{IP2}}$ are much smaller than $1$ for typical values $10 \leq \Fres \leq 1000$. 
% Theorem \ref{thm:IP2StabRes} holds 
% in any dimension $m$, but in the proof we will work out constants only for the 
% physically relevant values of $m$. 
Moreover, we emphasize that the stated results for \cref{ip1,ip2} are both for reconstructions from a \emph{single} diffraction pattern. Image recovery from \emph{two} holograms recorded at different distances  is treated in \secref{S5} as a corollary of the stability analysis of \cref{ip2}.

\vspace{1em}

\section{Stability analysis of \cref{ip1}} \label{S3}

\subsection{Principal approach} \label{SS3.1}

We start our analysis with \cref{ip1}, corresponding to the recovery of general complex-valued images from a single hologram without homogeneity constraint.  In order to understand its mathematical structure, it is instructive to rewrite the forward operator $T$ in the form
\begin{equation}
\begin{aligned}
 T h &= 2 \Re \big( \cD( h )  \big) = \cD(h) + \closure{\cD(h)} = \cD(h) + \cF^{-1} \big( \mF^{-1} \cdot  \cF(\closure h )  \big)  \\
     &= \cD(h) + \cD^{-1} ( \closure h)   \MTEXT{for}  h \in L^2(\mR^m).
     \end{aligned}\label{eq:3.1-1}
\end{equation}
 Here, the overbar denotes complex conjugation and we have used that the Fresnel propagation factor $\mF$ is unitary. According to \eqref{eq:3.1-1}, the linearized contrast $Th$ in the intensity data is given by a superposition of the propagated image $\cD(h)$ and the \emph{back-propagated twin-image} $\cD^{-1}(\closure h)$. Since $\cD$ is unitary, reconstructing $h$ from one of these components \emph{alone} would be straightforward. Solving  \cref{ip1} thus amounts to disentangling image and twin-image.

\begin{figure}
     \centering
     \includegraphics[width=\textwidth]{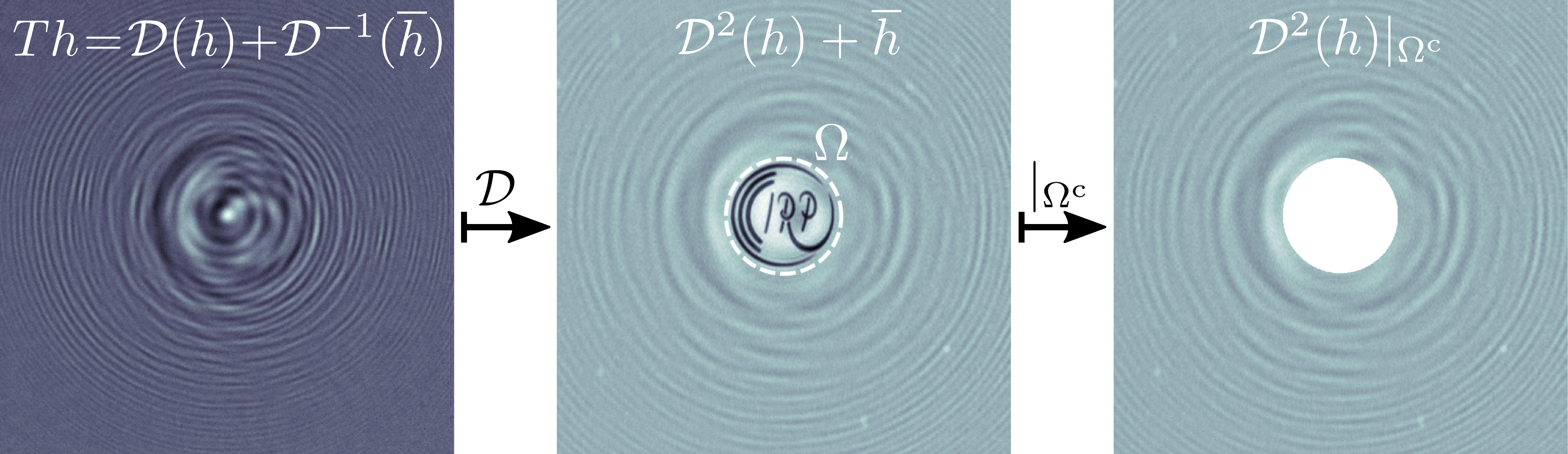}
     \caption{Illustration of the principal idea for the stability analysis of  \cref{ip1}. By applying the Fresnel propagator $\cD$ to data $Th$, the \emph{twin-image} $\closure h$ becomes sharp (\emph{Gabor holography}, see logo in central panel). By restricting to the complement $\Omega^{\textup c}$ of the support-domain $ \Omega \supset \supp(h)$, $\closure h$ is eliminated and incomplete Fresnel data $\cD^2(h)|_{\Omega^{\textup c}}$ is obtained (right panel). Images show real parts of numerically computed fields from a hologram (left panel) recorded at GINIX \cite{Salditt2015GINIX}, P10-beamline, DESY. \label{figure2}}
  \end{figure}

In order to separate these components, we propagate the near-field hologram $Th$ by application of $\cD$. By \eqref{eq:3.1-1}, this yields
% :
% By applying the Fresnel propagator $\cD$ to \eqref{eq:3.1-1}, we obtain
\begin{equation}
 \cD T h  = \cD^2( h) +  \closure h.  %\MTEXT{for all}  h \in L^2(\mR^m)
  \label{eq:3.1-2}
\end{equation}
Accordingly, we recover the sharp twin-image $\closure h$  up to perturbations originating from the doubly propagated image $\cD^2(h)$. This is the  principle of \emph{Gabor holography} \cite{Gabor1948_Holography}, which can be used as a \emph{qualitative} image reconstruction technique \cite{PaganinXRay} as illustrated in \cref{figure2}. Here, we follow this approach in a converse manner: rather than contenting ourselves with the perturbed twin-image $\closure h$, we exploit a \emph{support constraint}
\begin{equation}
 \supp(h) = \closure{\{ \bx \in \mR^m: h(\bx) \neq 0 \}} \subset \Omega
\end{equation}
by restricting \eqref{eq:3.1-2} to the complement of $\Omega$. This yields
\begin{equation}
 \cD T h|_{\comp \Omega}  = \cD^2( h)|_{\comp \Omega} + \closure h |_{\comp \Omega}
= \cD^2( h)|_{\comp \Omega} 
 \MTEXT{for any}  h \in L^2_\Omega. \label{eq:3.1-3}
\end{equation}

By the proposed propagation-and-restriction procedure, the twin-image is thus completely eliminated from the data as sketched in \cref{figure2}. Note that the map $Th \mapsto \cD(Th)|_{\comp \Omega} = \cD^2( h)|_{\comp \Omega}$ is norm-decreasing in $L^2(\mR^m)$. Hence, \eqref{eq:3.1-3} implies that the solution of \cref{ip1} is at most as ill-posed as the reconstruction from incomplete Fresnel data $\cD^2( h)|_{\comp \Omega}$. This reduction to a data completion problem is the principal idea of our stability analysis for Inverse Problem \ref{ip1}.

\vspace{.5em}
\subsection{Reduction to Fourier completion problem} \label{SS3.2}

In order to gain a simpler expression for the reduced data $\cD^2( h)|_{\comp \Omega}$ obtained in \secref{SS3.1}, we use an alternate form of the Fresnel propagator $\cD$. By application of the convolution theorem to \eqref{eq:FresnelProp}, the following representation can be obtained (see e.g.\ \cite{PaganinXRay}):
\begin{equation}
  \cD(h)(\bx)  = \E^{-\I m \pi /4} \tNF ^{\frac m 2 } \nF(\bx) \cdot \cF (\nF \cdot h)(\tNF \bx ), \;\; \nF(\bx) := \exp\left( \frac{\I}{2} \tNF |\bx|^2 \right) \label{eq:FresnelPropAltForm}
\end{equation}
for all $\bx \in \mR^m$. 
%Here and in the following, we employ the \emph{modified Fresnel number} $\tNF := 2\pi \Fres$ rather than the more common parameter $\Fres$ in order to avoid excessive occurence of $2\pi$-factors. 
\eqref{eq:FresnelPropAltForm} reveals that - up to pointwise multiplications with the unitary factor $\nF$ and rescaling - $\cD$ may be written as a Fourier transform. Combined with the approach outlined in \secref{SS3.1}, this allows to identify stability of \cref{ip1} with the reconstruction of a function from \emph{incomplete Fourier data}:

\vspace{.5em}
\begin{theorem}[Reduction to Fourier data completion problem] \label{thm:IP1FourCompRed}
Within the setting of \cref{ip1}, let $\OmegaF := \{ (\tNF/2) \bx : \bx \in \Omega\}$ for arbitrary $\Omega \subset \mR^m$. Then
\begin{equation}
 \norm{ T h } \geq \bnorm{ \cF ( \nF^{1/2} \cdot h )|_{ \comp \OmegaF} } \MTEXT{for all} h \in L^2_\Omega \label{eq:IP1FourCompRed-1}
\end{equation}
where $\nF^{1/2}(\bx) =  \exp\left( \I \tNF |\bx|^2 / 4\right)$. In particular, we have the relative stability estimate
\begin{equation}
 C_{\textup{IP1}}(\Omega,\tNF) = \inf_{h\in L^2_\Omega, \norm{h} = 1 } \norm{Th} \geq \inf_{h\in L^2_\Omega, \norm{h} = 1 }  \norm{ \cF ( h )|_{ \comp \OmegaF} } \label{eq:IP1FourCompRed-2}
\end{equation}
\end{theorem}
\vspace{.5em}
\begin{proof}
$\cD$ is unitary and  the restriction to $\comp \Omega \subset \mR^m$ defines an orthogonal projection in $L^2(\mR^m)$. Hence, \eqref{eq:3.1-3} implies the estimate
\begin{equation}
 \norm{Th } = \norm{ \cD T h } \geq \norm{ \cD T h|_{\comp \Omega} } = \norm{\cD^2( h)|_{\comp \Omega}}  \MTEXT{for any} h \in L^2_\Omega. \label{eq:IP1FourCompRed-proof-1}
\end{equation} 
By \eqref{eq:FresnelProp}, $\cD^2$ is again a Fresnel propagator, yet to the Fresnel number $\tNF/2$. Accordingly, employing the alternative form \eqref{eq:FresnelPropAltForm} and exploiting that $|\E^{-\I m \pi /4}\nF^{1/2}|\equiv 1$ gives
\begin{equation}
 \norm{\cD^2( h)|_{\comp \Omega}} = \left( \tfrac{1}{2} \tNF \right)^{\frac m 2 } \norm{ \cF (\nF^{1/2} \cdot h)\big( \tfrac{1}{2} \tNF \cdot \big)|_{\comp \Omega} }. \label{eq:IP1FourCompRed-proof-2}
\end{equation}
Introducing new coordinates $\bxi := (\tNF/2) \bx$   and using that $\comp \OmegaF = (\tNF/2) \cdot \comp \Omega$ holds by definition, this expression can be simplified to
\begin{equation}
\begin{aligned}
 \norm{\cD^2( h)|_{\comp \Omega}}^2 &=\left( \tfrac{1}{2} \tNF \right)^{m} \int_{\comp \Omega} \big|\cF (\nF^{1/2} \cdot h)( \tfrac{1}{2} \tNF \bx)|^2 \; \D\bx  \\
&= \int_{\comp \OmegaF} \big|\cF (\nF^{1/2} \cdot h) (\bxi)|^2 \; \D\bxi = \norm{\cF( \nF^{1/2} \cdot h )|_{\comp \OmegaF}}^2. \label{eq:IP1FourCompRed-proof-3}
\end{aligned}
\end{equation}
Combining \eqref{eq:IP1FourCompRed-proof-1} and \eqref{eq:IP1FourCompRed-proof-3} yields the first assertion  \eqref{eq:IP1FourCompRed-1}. The second estimate \eqref{eq:IP1FourCompRed-2} then follows from the fact that the map $h \mapsto \nF^{1/2} \cdot h$ is isometric and bijective on $L^2_\Omega$.
\end{proof}
\vspace{.5em}

\cref{thm:IP1FourCompRed} states that the solution of \cref{ip1} is \emph{at least as stable} as the reconstruction of an $L^2$-function $f$ with support in ${\Omega} \subset \mR^m$ from incomplete Fourier measurements $\cF(f)|_{{\comp \OmegaF}}$. For compact $\Omega$, the latter Fourier completion problem can be shown to be unique and even \emph{well-posed} by employing an uncertainty principle for the 1D-Fourier transform derived by Nazarov \cite{Nazarov1993CompactFTBound} (see \cite{HavinJoericke1994UncertPrinc} for an English proof and \cite{Jaming2007_NazarovND} for its multidimensional generalization). Rather than following this approach for general support shapes, however, we will restrict to the special case of rectangular support-domains $\Omega$. This will enable a more explicit characterization of the dependence of $C_{\textup{IP1}}$ on the Fresnel number $\tNF$ via \cref{thm:IP1FourCompRed} as well as additional insights concerning the nature of the least stable modes of the forward operator $T$.

 \vspace{.5em}
 \subsection{Stability result for stripe-shaped supports} \label{SS3.3}

 In the following, we restrict to the simple case of a stripe-shaped support-domain $\Omega \subset \mR^m$ as considered in \cref{thm:IP1StabRes}. Note that the forward operator $T$ is rotationally- and translationally invariant as is the Fresnel propagator $\cD$. Hence, it is sufficient to consider domains of the form $\Omega := [-b/2 ; b/2] \times \mR^{m-1}$. Moreover, $b$ can be set to $1$ which means that we define the Fresnel number $\tNF$ with respect to the support diameter. Accordingly, we may indeed restrict to the special case
 $ \Omega := [-1/2 ; 1/2] \times \mR^{m-1}$
 as done in \cref{thm:IP1StabRes} without loss of generality.

 Now we can employ the characterization of the stability constant in \cref{thm:IP1FourCompRed} for this special domain. We define $I := [-1/2 ; 1/2]$ and $\IF := [-\tNF/4 ; \tNF/4]$ for notational convenience.
Recall that the $m$-dimensional Fourier transform $\cF^{(m)}$ is a tensor 
product of one-dimensional Fourier transforms $\cF^{(1)}$ applied along the different coordinate dimensions, i.e.\ 
$\cF^{(m)}=\cF^{(1)}\otimes \cdots\otimes\cF^{(1)}$.
Owing to the cartesian product structure of $\Omega = I \times \mR^{m-1}$ and $\comp \OmegaF = \comp{\IF} \times \mR^{m-1}$, the restricted Fourier transforms 
\begin{equation}
\cF^{(m)}_{\Omega,\comp\OmegaF}\,:\,L^2_{\Omega}\to L^2_{\comp\OmegaF};\; h \mapsto (\cF^{(m)}h)|_{\comp\OmegaF},\qquad 
\cF^{(1)}_{I,\comp\IF}\,:\, L^2_I\to L^2_{\comp\IF};\; h \mapsto (\cF^{(1)}h)|_{\comp\IF},   \label{eq:IP1DefFF}
\end{equation}
are likewise related by $\cF^{(m)}_{\Omega,\comp\OmegaF} = \cF^{(1)}_{I,\comp{\IF}}\otimes \cF^{(1)}\otimes \cdots  \otimes \cF^{(1)}$. Applying this relation to the stability estimate in \eqref{eq:IP1FourCompRed-2} and exploiting unitarity of the Fourier transform yields
\begin{eqnarray}
C_{\textup{IP1}}(\Omega,\tNF) 
&=& \inf_{h\in L^2_\Omega, \norm{h} = 1 }   \bnorm{   \cF^{(m)}_{\Omega,\comp \OmegaF} h}
= \inf_{h\in L^2_I, \norm{h} = 1 }   \bnorm{\cF^{(1)}_{I,\comp{\IF} }h }\nonumber\\
 &=& 
\bigg(1 - \sup_{h\in L^2_I, \norm{h} = 1 }   \big \|   \cF^{(1)}  ( h )  |_{ \IF  } \big\|^2 \bigg)^{1/2} . \label{eq:CIP1_PrincSingVal}
 \end{eqnarray}
 We thus need to estimate the norm of the 1D-operator $\FF: L^2_I \to L^2_{\IF};\, h \mapsto \cF^{(1)}(h)|_{\IF}$. %:= \cF^{(1)}_{I, \IF}$.
 This is achieved by explicit computation of the operator $\FF^\ast \FF$ ($\FF^\ast$: adjoint of $\FF$), which turns out to be part of a well-studied family of compact and self-adjoint integral operators. 
Their eigenfunctions are known as prolate spheroidal wave 
functions, and the eigenpairs have been studied for example in \cite{Slepian1983FourierCompactSpectrum,Slepian1965_CpctFTSVD}.
By applying these known results we obtain the following theorem:

  \vspace{.5em}
  \begin{theorem}[Spectral characterization of $\FF^\ast \FF$] \label{thm:SVDFF}
 Let $\FF$ defined by \textup{\eqref{eq:IP1DefFF}}. Then 
   \begin{eqnarray}
 \FF^\ast \FF(f)\left(\frac x 2 \right) =   \int_{1}^{1}  \frac{ \sin \left( c (x-y) \right) }{ \pi (x-y) } f\left(\frac y 2 \right) \; \D y \MTEXT{with} c := \tNF / 8. \label{eq:IP1ProlateOP}
 \end{eqnarray}
for all $h \in L^2_I$ and $x \in [-1;1]$, and $\FF^\ast \FF$ is compact. 
The eigenvalues $\{ \lambda_{c,j} \}_{j \in \mN_0} \subset \mR_+$ and associated eigenfunctions $\{ \psi_{c,j} \}_{j \in \mN_0} \subset L^2_I$ of $\FF^\ast \FF$ thus coincide with those in \textup{\cite{Slepian1983FourierCompactSpectrum,Slepian1965_CpctFTSVD}}. In particular, all eigenvalues $\lambda_{c,0} > \lambda_{c,1} > \ldots$ have multiplicity one, and the $\psi_{c,j}$ may be chosen to form an orthonormal basis of $L^2_I$. Moreover, $\lambda_{c,0} < 1$ holds true and, for fixed $j \in \mN_0$ and $\tNF\to \infty$, $\lambda_{c,j}$ has the asymptotic expansion
 \begin{equation}
  1 - \lambda_{c,j}  = \frac{ (2\pi)^{\frac 1 2} \tNF ^{j + \frac 1 2 }}{j!} \left( 1 - \frac{6j^2 -2j + 3 }{4 \tNF}  + \cO \left( \frac 1 {\tNF^2 } \right)\right) \exp(-\tNF /4 ).  \label{eq:IP1EVProlateAsymp-1}
 \end{equation}
 \end{theorem}
  \vspace{.5em}
 \begin{proof}
  The restriction to the interval $\IF = [-\tNF/4 ; \tNF/4]$ can be written in the form of a multiplication with its indicator function $\boldsymbol 1 _ {\IF}$. By the convolution theorem, we thus obtain for all $h \in  L^2_I  $, $x \in [-1; 1]$
  \vspace{.5em}
  \begin{eqnarray*}
 \FF^\ast \FF (h)\left(\frac x 2 \right) &=& \cF^{-1} \left( \boldsymbol 1 _ {\IF} \cdot  \cF( h )  \right)  \left(\frac x 2 \right) = (2\pi)^{- \frac 1 2} \cF^{-1}\big(\boldsymbol 1 _ {\IF} \big) \ast h  \left(\frac x 2 \right) \nonumber \\
 &=& \frac 1 { 2 \pi } \int_{-1/2}^{1/2}  \frac{ 2 \sin \big( \frac{\tNF} 4 (\frac x 2 -y) \big) }{\frac 1 2 x-y   } f(y) \; \D y 
 =   \int_{1}^{1}  \frac{ \sin \left( c (x-y) \right) }{ \pi (x-y) } f\left(\frac y 2 \right) \; \D y. 
 \end{eqnarray*}
The spectral characterization of the resulting integral operators in \cite{Slepian1983FourierCompactSpectrum,Slepian1965_CpctFTSVD} directly yields the claimed properties of the eigensystem $\{ (\lambda_{c,j}, \psi_{c,j}) \}_{j \in \mN_0}$  of  $\FF^\ast \FF$. In particular, the asymptotic expansion \eqref{eq:IP1EVProlateAsymp-1} is an analogue of the formula \cite[eq.\ (2)]{Slepian1983FourierCompactSpectrum}.
 
 Since $\FF$ is a restriction of the Fourier transform, we  have $\norm{\FF} \leq \norm{\cF}  =1$. Hence, the principal eigenvalue of $\FF^\ast \FF$ must satisfy $\lambda_{c,0} \leq 1$. If $\lambda_{c,0} =1$ then
 \begin{equation*}
  \norm{\cF(\psi_{c,0})|_{\comp{\IF} } }^2 
	= \norm{\cF(\psi_{c,0})  }^2  -  \norm{\cF(\psi_{c,0})|_{ \IF  } }^2
	= \norm{\psi_{c,0}}^2- \lambda_{c,0} \norm{\psi_{c,0}}^2 
	= 0,
 \end{equation*}
 i.e.\ $\cF(\psi_{c,0})$ would have to vanish outside the interval $\IF$. However, as $\psi_{c,0} $ is compactly supported, $\cF(\psi_{c,0})$ is an entire function and thus vanishes identically if $\cF(\psi_{c,0})|_{\comp{\IF}} = 0$. This is impossible since $\psi_{c,0} $ is an eigenfunction. Hence, $\lambda_{c,0} < 1$ must hold true.
 \end{proof}
   \vspace{.5em}

We emphasize the nontrivial dependence of both the eigenvalues $\lambda_{c,j}$ and the eigenfunctions $\psi_{c,j}$  on the parameter $c = \tNF/8$. For convenience, however, we will suppress the subscript $c$ in the following.
 \cref{thm:SVDFF} constitutes the final ingredient which is needed to prove the sought stability result for \cref{ip1}:

 \vspace{.5em}
 \begin{proof}[Proof of \cref{thm:IP1StabRes}] 
 According to the characterization of the stability constant in \eqref{eq:CIP1_PrincSingVal}, $C_{\textup{IP1}}(\Omega,\tNF) $ can be expressed in terms of the operator norm of $\FF$. Since $\FF^\ast \FF$ is compact with principal eigenvalue $\lambda_0 < 1 $ and orthonormal eigenfunctions $\{\psi_j\}_{j\in \mN_0}$ as characterized in \cref{thm:SVDFF}, we have
 \begin{equation*}
  \norm{\FF}^2 = \norm{\FF^\ast \FF} = \sup_{j \in \mN_0} \norm{ \FF^\ast \FF \psi_j }  = \sup_{j \in \mN_0} \lambda_j = \lambda_0
 \end{equation*}
 By \eqref{eq:CIP1_PrincSingVal}, this implies $C_{\textup{IP1}}(\Omega,\tNF)^2 = 1 - \lambda_0 > 0$ , i.e.\ \emph{well-posedness} of \cref{ip1}. Setting $j = 0$ in \eqref{eq:IP1EVProlateAsymp-1} yields the asymptotic characterization \eqref{eq:IP1StabRes-2}.
 \end{proof}
 \vspace{.5em}

 \subsection{Characterization of the least stable modes} \label{SS3.4}

So far, we have not exploited the full potential of the reduction to a Fourier data completion problem in \cref{thm:IP1FourCompRed} yet: only the characterization of the stability constant \eqref{eq:IP1FourCompRed-2}, i.e.\ of the \emph{worst-case-stability}, has been used in the proof of \cref{thm:IP1StabRes}. Notably however,  the more general estimate \eqref{eq:IP1FourCompRed-1} even bounds the contrast $\norm{Th}$ attained by \emph{individual images} $h$ with respect to the corresponding incomplete Fourier data $ \cF(\nF^{1/2} \cdot h) | _{\comp \OmegaF } $. This enables a precise prediction of the reconstruction stability for different image modes in  \cref{ip1}  beyond the universal lower bound proven in \cref{thm:IP1StabRes}.

In order to avoid notational difficulties in the argument, we do the analysis for a box-shaped support-domain $\Omega := [-1/2; 1/2]^m = I^m$. Owing to the simple Cartesian geometry, the results obtained for a stripe support are easily generalized to this case, including a characterization of the stability of individual modes. 
We define $m$-dimensional prolate spheroidal wave functions as the tensor product 
\begin{equation} 
 \psi_{\bj}(\bx) := (\psi_{j_1} \otimes \ldots \otimes \psi_{j_m})(\bx)
:= \prod_{j = 1}^m \psi_j(x_j) 
\MTEXT{for} \bj = (j_1, \ldots, j_m) \in \mN_0^m.  \label{eq:IP1MDProlate}
\end{equation}
Moreover, let $\ip{f}{g} := \int_\Omega f(\bx) \closure{ g(\bx)}\; \D \bx$ for $h_1, h_2\in L^2_\Omega$ denote usual $L^2$-inner product. With this notation, we obtain the following modal stability estimates:

\vspace{.5em}
  \begin{theorem}[Stability of individual modes in \cref{ip1}] \label{thm:LeastStableModes}
 Let $ \Omega = I^m$ and let $\{(\lambda_j,\psi_j)\}_{j\in \mN_0}$ denote the eigenvalue decomposition of $\FF^\ast \FF$ in \cref{thm:SVDFF}. 
 Moreover, define $\phi_{\bj}   := \nF^{-1/2} \cdot    \psi_{\bj}$ for all multi-indices $ \bj \in \mN_0^m$.
 Then $\{\phi_{\bj}\}_{\bj \in \mN_0^m}$ is an orthonormal basis of $L^2_\Omega$, 
and  with  $c_{\textup{IP1}, \bj} :=  (1- \prod_{l=1}^m \lambda_{j_l})^{1/2}$ we have
 \begin{equation}
 \norm{Th}^2 \geq \sum_{\bj \in \mN_0^m} c_{\textup{IP1}, \bj}^2  |\ip{h}{\phi_{\bj}}|^2    \MTEXT{for any} h = \sum_{\bj \in \mN_0^m} \ip{h}{\phi_{\bj}} \phi_{\bj} \in L^2_\Omega. \label{eq:LeastStableModes-1}
 \end{equation}
 \end{theorem}
  \vspace{.5em}
\begin{proof}
  As  $\{\psi_{j}\}_{j\in \mN_0}$ is an orthonormal basis of $L^2_I$, the tensor products $\psi_{\bj} = \psi_{j_1} \otimes \ldots \otimes \psi_{j_m} $ form an orthonormal basis of $L^2_\Omega$ with $\Omega = I^m$. Since $\nF^{-1/2}$ is unitary, the same is true for $\{\phi_{\bj}\}_{\bj \in \mN_0^m}$ so that  any $h \in L^2_\Omega$ can be written as $ h = \sum_{\bj \in \mN_0^m} a_{\bj} \phi_{\bj}$ with $a_{\bj} := \ip{h}{\phi_{\bj}}$.
  By the tensor product structure of $\psi_{\bj}$ and $\cF^{(m)}$, we further have
  \begin{align*}
   \cF^{(m)} ( \nF^{1/2} \cdot \phi_{\bj} )|_{\OmegaF} &= \cF^{(m)}( \psi_{\bj} )|_{\IF^m}  
	%=  \cF^{(1)}  \ldots    \cF^{(m)} ( \psi_{j_1} \otimes \ldots \otimes  \psi_{j_m} )|_{\IF^m} \nonumber \\
   = (\cF^{(1)}  \psi_{j_1} )|_{\IF} \otimes \ldots \otimes (\cF^{(1)}  \psi_{j_m} )|_{\IF} \\
   &= \FF\psi_{j_1} \otimes \ldots \otimes \FF\psi_{j_m}  
%	\label{eq:LeastStableModes-proof1}
  \end{align*}
  for all $\bj = (j_1, \ldots, j_m) \in \mN_0$. As the $\psi_{l}$ are eigenfunctions of the 1D-map $ \FF^\ast \FF$ to the eigenvalues $\lambda_l$ (see \cref{thm:SVDFF}), the above relation implies
  \begin{eqnarray}
    &&\left\langle \cF ( \nF^{1/2} \cdot \phi_{\bj})|_{  \OmegaF}, \cF ( \nF^{1/2} \cdot \phi_{\bk} )|_{\OmegaF} \right\rangle \label{eq:LeastStableModes-proof2}\\
		&=&  \prod_{l=1}^m \left\langle  \FF^\ast \FF(\psi_{j_l} ),   
    \psi_{k_l}  \right\rangle 
    = \prod_{l=1}^m \lambda_{j_l} \left\langle  \psi_{j_l},  \psi_{k_l}\right\rangle 
		=  \delta_{\bj\bk} \prod_{l=1}^m \lambda_{j_l} 
    = \delta_{\bj\bk} ( 1- c_{\textup{IP1}, \bj}^2)  \nonumber% \cases{\prod_{l=1}^m \sigma_{k_l}^2 &\text{if } \bj =\bk \\ 0 & \text{else}}
  \end{eqnarray}
  for all $ \bj, \bk \in \mN_0^m$ where $\delta_{\bj\bk} \in \{0,1\}$ is the multidimensional Kronecker-symbol.
  Using the principal bound \eqref{eq:IP1FourCompRed-1} from \cref{thm:IP1FourCompRed}, unitarity of  $\cF$ and \eqref{eq:LeastStableModes-proof2} finally gives
\begin{align*}
 \norm{Th}^2 &\geq \norm{ \cF ( \nF^{1/2} \cdot h )|_{\comp \OmegaF} }^2  = \norm{h}^2 -  \norm{ \cF ( \nF^{1/2} \cdot h)|_{ \OmegaF} }^2 \\
 &= \sum_{\bj \in \mN_0^m}  |a_{\bj}|^2  \!-\!\!\! \sum_{\bj, \bk \in \mN_0^m} \! a_{\bj} \cc{a_{\bk} } \bigg\langle \cF ( \nF^{1/2} \!\cdot\! \phi_{\bj})|_{ \OmegaF}, \cF ( \nF^{1/2} \!\cdot\! \phi_{\bk} )|_{\OmegaF} \bigg\rangle 
 =  \sum_{\bj \in \mN_0^m}  c_{\textup{IP1}, \bj}^2 |a_{\bj}|^2. 
\end{align*} 
\end{proof}
\vspace{.5em}

In \cref{thm:LeastStableModes} we obtain \emph{individual} stability constants $ c_{\textup{IP1}, \bj}$ such that each of the constructed orthonormal basis modes $\phi_{\bj}$ attains data contrast $\norm{T \phi_{\bj} } \geq c_{\textup{IP1}, \bj}$. As the sequence of eigenvalues $\{\lambda_j\}_{k \in \mN_0}$ in \cref{thm:SVDFF} is strictly decreasing, it is readily seen that stability \emph{increases} with the multi-index $\bj$, i.e.\
\begin{equation}
 c_{\textup{IP1}, \bk} > c_{\textup{IP1}, \bj} \MTEXT{if} \bk \neq \bj, \; k_l \geq j_l \MTEXT{for all} l = 1,\ldots,m.
\end{equation}
According to \cref{thm:LeastStableModes}, the least stable modes  are thus exactly the basis functions  $\phi_{\bj}$ of (componentwise) small multi-index $\bj$, which are given by prolate spheroidal wave functions $\psi_ {\bj} = \psi_{j_1} \otimes \ldots \otimes \psi_{j_m}$, modulated by the Fresnel factor $\nF^{-1/2}$. Details on the shape of the $\psi_ {j}$ in turn are readily available, see e.g.\ \cite{Slepian1983FourierCompactSpectrum}. In particular, the index $j$ can be interpreted as a \emph{frequency} since $\psi_j$ is smooth and real-valued with exactly $j$ zeros and $j+1$ extrema within the interval $(-1 /2 ; 1/2)$. By extending this observation to the $\phi_{\bj}$ we deduce the nature of the least stable modes in  \cref{ip1}:

\vspace{.5em}
\begin{corollary}[Least stable modes in \cref{ip1}]\label{cor:LeastStableModes}
The least stable modes in \cref{ip1} for  $\Omega = [-1 /2; 1/2]^m$ 
are \emph{low-frequency} prolate spheroidal wave functions modulated by the Fresnel-factor $\nF^{-1/2}(\bx) = \exp(-\I\tNF|\bx|^2/4)$.
\end{corollary}
\vspace{.5em}

\subsection{Numerical Validation} \label{SS3.5}
 
 We return to the starting point of our stability analysis, namely the reduction to a Fourier data completion problem via the principal estimate \eqref{eq:IP1FourCompRed-1}. The question whether or not the derived stability bounds in  \cref{thm:IP1StabRes} and \cref{thm:LeastStableModes} are \emph{optimal} (or at least close to) crucially depends on the sharpness of this inequality. In the following we investigate this remaining issue numerically.
 
\begin{figure}
     \centering
     \includegraphics[width=\textwidth]{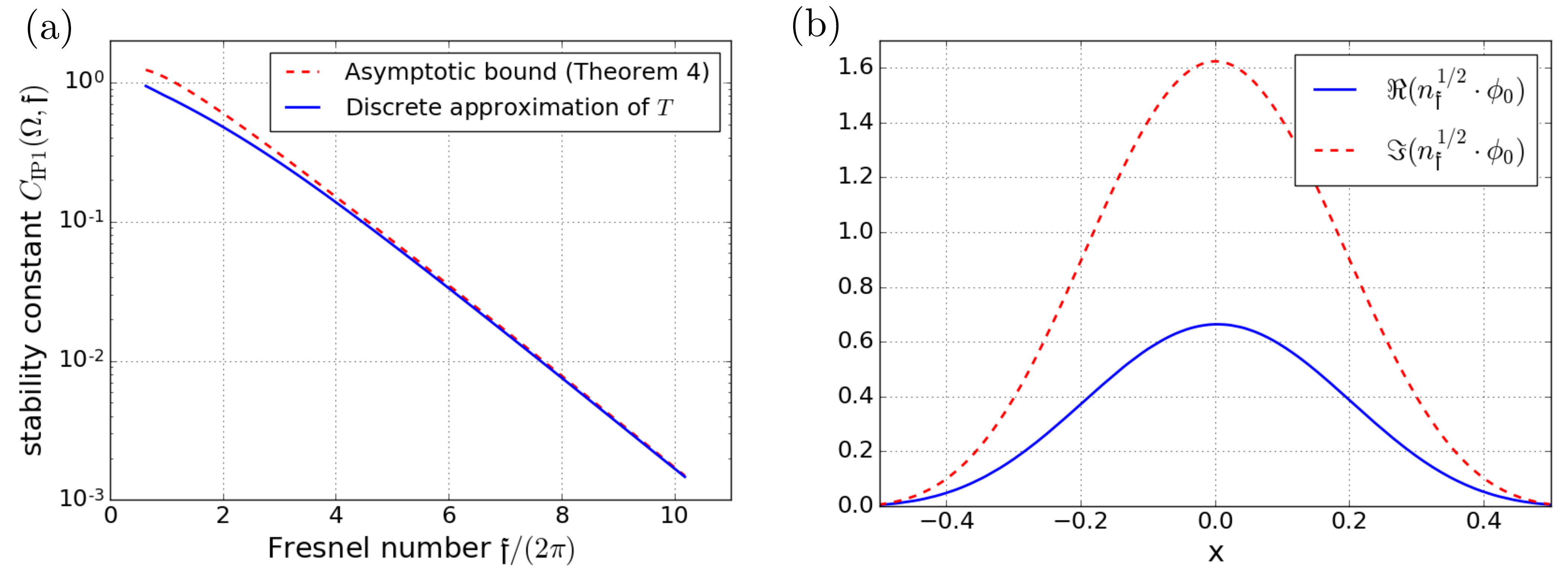}
     \caption{Numerical validation of the stability analysis for \cref{ip1}. \\ 
     (a) Comparison of the analytical bound \eqref{eq:IP1StabRes-2} for the stability constant $C_{\textup{IP1}}$ (red-dashed line) to numerical computations of the smallest singular value of $T$ (blue-solid). \\ % (discretized by FFTs on $512^2$ grid points) \\
     (b) Numerically computed least stable singular mode $\phi_0$ for $\Fres = \tNF/(2\pi) = 10$. The plotted modulation with $\nF^{1/2}$ reveals a unimodal structure as predicted by \cref{thm:LeastStableModes}. 
     \label{figure4}}
  \end{figure}

 To this end, we approximate the stability bound $C_{\textup{IP1}}( \Omega,\tNF)$ by computing the smallest singular value of a discretized forward operator $T$ in $m = 1$ dimensions. The Fresnel propagator is approximated by fast Fourier transforms on a 1D-grid of $512^2$ equidistant points, where the central 512 grid points form the supporting interval $\Omega = [-1/2; 1/2]$. We compute the smallest singular value of $T$ via a power method for Fresnel numbers $\Fres \in [1; 10]$. These numerical results for the complete, yet discretized forward operator, are compared to the asymptotic stability bound \eqref{eq:IP1StabRes-2} in \cref{thm:IP1StabRes}, neglecting the $\cO(\tNF^{-2})$-contributions. The different predictions for the stability constant $C_{\textup{IP1}}( \Omega,\tNF)$ are plotted in \cref{figure4}(a).
  
  The semilogarithmic plot shows excellent agreement between the analytical bound \eqref{eq:IP1StabRes-2} and the numerical approximation in the asymptotic limit $\tNF \to \infty$. This indicates that our stability analysis, based on the potentially lossy estimate \eqref{eq:IP1FourCompRed-1}, is surprisingly sharp. % - despite the discarding of ``half of the data'' by the underlying propagate-and-restrict approach in \cref{figure2}.
  Accordingly, one might expect that also the corresponding least stable modes $\phi_{\bj}$ are well-characterized by \cref{thm:LeastStableModes}. This is supported by the simulation results: \cref{figure4}(b) exemplarily plots the numerically computed  mode $\phi_0$ to the minimum singular value of $T$ for $\Fres = \tNF/(2\pi) = 10$. According to \cref{thm:IP1StabRes}, the plotted product with the factor $\nF^{1/2}$ should yield the zeroth order prolate spheroidal wave function $\psi_0$. This is confirmed by the smooth unimodal profiles obtained in \cref{figure4}(b).

\vspace{1em}

\section{Stability analysis of \cref{ip2}} \label{S4}

%As the forward operator $S_{\alpha}$ in \eqref{eq:fwModel2} arises from $T$ by a mere \emph{restriction} to homogeneous images, well-posedness and stability of \cref{ip1} directly carry over to \cref{ip2}. However, it should be emphasized that the proven \emph{exponential} rates for the constant $C_{\textup{IP1}}(\Omega,\tNF) \gtrsim \exp(-\tNF / 8)$ hardly guarantee stability in a practical sense for Fresnel numbers $\Fres \geq 100$, as encountered in typical near-field imaging experiments. We have seen in \secref{SS:StabResults} that \cref{ip2} is in general much less ill-posed than \cref{ip1}. In this spirit, \cref{thm:IP2StabRes} establishes improved \emph{algebraic} rates for the stability constant $C_{\textup{IP2}}(\Omega,\tNF) \gtrsim \tNF^{-1}$ with the Fresnel number $\tNF$. The aim of the following analysis is to prove this stability result.

The aim of this section is to prove \cref{thm:IP2StabRes}, establishing 
\emph{algebraic} rates of the stability constant $C_{\textup{IP2}}(\Omega,\tNF) \gtrsim \tNF^{-1}$ under a  homogeneity constraint for the imaged object. % for homogeneous objects  with the Fresnel number $\tNF$. 

\vspace{.5em}
\subsection{Preparations and Fourier domain splitting} \label{SS4.1}

Throughout this section, let $\varphi \in L^2(\mR^m)$ and let $\hat \varphi := \cF(\varphi)$ denote its Fourier transform. Recall from definition \eqref{eq:fwModel2} that 
the operator $S_{\alpha}$ denotes multiplication with  
$2s_\alpha(\bxi) =  2\sin \left( |\bxi|^2/(2\tNF) + \alpha \right)$ in the 
Fourier domain, and hence 
\begin{equation}
 \norm{S_{\alpha}\varphi}  =  \norm{\cF S_{\alpha}\varphi} = 2 \norm{ s_\alpha \cdot \hat \varphi  }. \label{eq:IP2Contrast}
\end{equation}
According to \eqref{eq:IP2Contrast}, the images $\varphi$ attaining low contrast, i.e.\ small  $\norm{S_{\alpha}\varphi}$, are exactly those for which $\hat \varphi$ is concentrated about the zero-manifolds of the CTF $s_\alpha$. As signals $\hat \varphi \in L^2(\mR^m)$ may be arbitrarily sharply peaked about these manifolds of zero contrast, \cref{ip2} is ill-posed for general images $\varphi \in  L^2(\mR^m)$.
  \begin{figure}
     \centering
     \includegraphics[width=\textwidth]{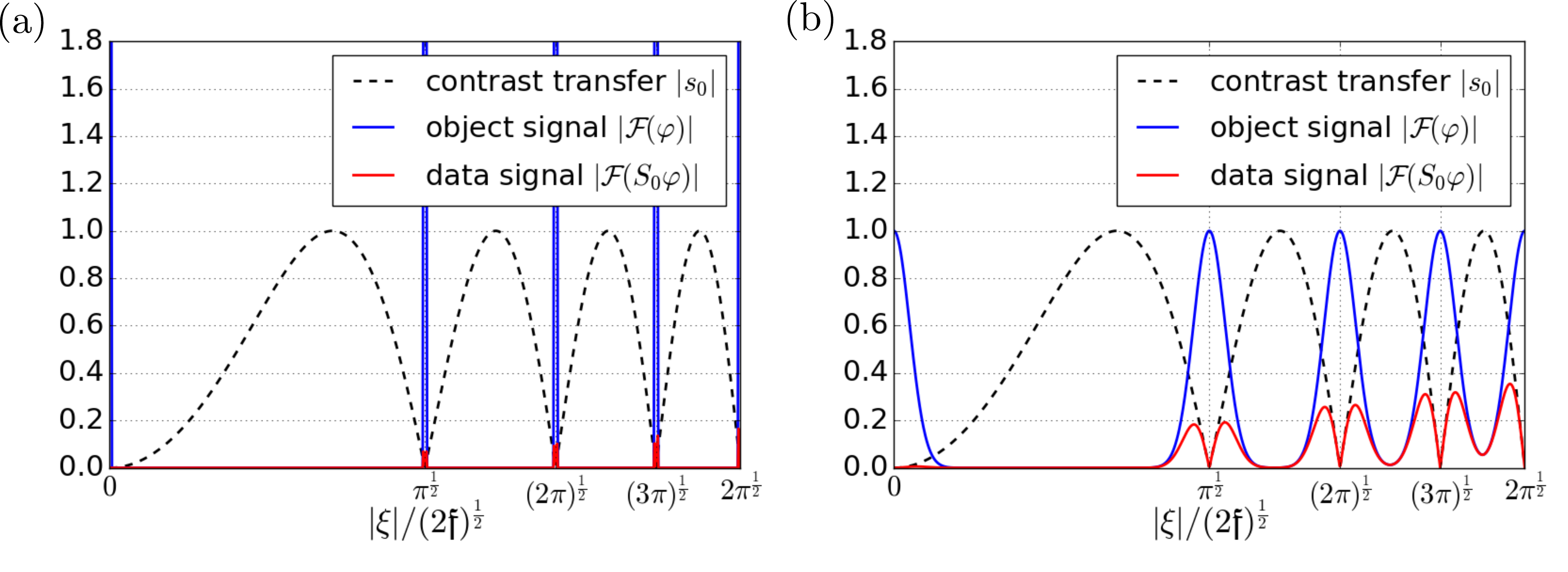}
     \caption{Illustration of the principal argument for stability of \cref{ip2}: the black-dashed line plots the radially symmetric contrast transfer function $s_\alpha$ for $\alpha = 0$. Blue and red solid lines show examples of Fourier space signals of objects $\varphi$ and their images under the forward operator $S_{\alpha}$, respectively. \\ 
     (a) General images $\varphi \in L^2(\mR^m)$: $\cF(\varphi)$ may be arbitrarily peaked at the zeros of $s_\alpha$ \\
     (b) Support constraint $\varphi \in L^2_\Omega$: $\cF(\varphi)$ is smooth and has finitely sharp peaks, which ensures minimum data contrast $\norm{\cF (S_{\alpha}\varphi)} \geq C_{\textup{IP2}} \norm{\cF (\varphi)}$ 
     \label{figure5}}
  \end{figure}
  
Now, we additionally assume a support constraint $\varphi \in L^2_\Omega$ for $\Omega := \closure { B(0,R) }$. Why does this constraint ensure well-posedness  in the light of the problematic CTF-zeros? The explanation lies in the well-known fact that a compact support in real-space implies $\sC^\infty$-smoothness (indeed analyticity) of the Fourier transform with norm-bounds on the derivatives in terms of the support size $R$. Owing to this regularity, $\hat \varphi$ may not be arbitrarily concentrated about the zero-manifolds of $s_\alpha$, which enables stability as illustrated in \cref{figure5}. 
The following lemma quantifies the smoothness of $\hat \varphi$ in a suitable form for the subsequent analysis:

\vspace{.5em}
\begin{lemma} \label{lem:CpctSuppDer2Bound}
Let $g \in L^2_\Omega$ with support in $\Omega = \overline{B(0, R)}$ and Fourier transform $\hat g := \cF( g)$. Let $\Delta$ be the Laplacian on $\mR^m$. Then we have for any measurable set $B \subset \mR^m$
\begin{equation}
 \int_{B} - \Delta |\hat g|^2 \; \D\bxi \leq 2 R^2 \norm{g} ^2.  \label{eq:CpctSuppDer2Bound}
\end{equation}
\end{lemma}
\vspace{.5em}
\begin{proof}
 The compact support of $g$ implies infinite smoothness of $\hat g = \cF( g)$.  
Using the identity $-\Delta \cF(g) = \cF(|\bx|^2 g)$ and Cauchy-Schwarz's inequality 
we obtain
 \begin{eqnarray*}
  \int_{B} - \Delta |\hat g|^2 \; \D\bxi &=& -2 \int_B \Big(|\nabla \hat g|^2 + \; \Re \big( \overline{\hat g} \cdot \Delta \hat g \big) \Big) \; \D\bxi  \leq 2 \left |  \int_B   \overline{\hat g} \cdot  \cF(|\bx|^2  g)  \; \D\bxi \right|  \\
  &\leq& 2\norm{\hat g|_B} \norm{\cF(|\bx|^2  g)|_B} \leq   2\norm{ \hat g} \norm{\cF(|\bx|^2  g)} \leq  2\norm{g}  \norm{|\bx|^2  g}.
\end{eqnarray*}
As $ g$ vanishes outside $\overline{B(0, R)}$, we further have $\norm{|\bx|^2 g}  \leq R ^2  \norm{ g}$.
\end{proof}
\vspace{.5em}

\noindent Note that Lemma \ref{lem:CpctSuppDer2Bound} can be interpreted as an \emph{uncertainty principle}: a bound for the derivative $-\Delta |\cF(g)|^2$, limiting the sharpness of features in Fourier space, arises from the confinement $\supp( g) \subset \overline{B(0, R)}$ of the corresponding real-space signal $g$.

For a quantitative analyis, we decompose the norm on the right-hand side of \eqref{eq:IP2Contrast} into stable bulk integrals and potentially unstable parts about the CTF-zeros by cutting out the subdomain $D_{\varepsilon} := \{ \bxi \in \mR^m : |s_\alpha(\bxi)| \geq \sin(\varepsilon) \}$ for some $ 0 < \varepsilon \leq \pi/6$:
\begin{equation}
\norm{ s_\alpha \cdot \hat \varphi }^2  =  \underbrace{ \int_{D_{\varepsilon} }  s_\alpha ( \bxi )^2  |\hat \varphi (\bxi)|^2  \; \D   \xi }_{ =: J_{\varepsilon} } + \sum_{j = 0}^\infty  \underbrace{ \int_{B_j}   s_\alpha ( \bxi )^2  |\hat \varphi(\bxi)|^2 \; \D   \xi   }_{=: J_j} \label{eq:IP2StabDecomp} 
\end{equation}
Here, the $B_j$ denote the annular connected components  of $ \mR^m \setminus D_{\varepsilon}$ about the $j$-th zero-manifold of $s_\alpha$ at radius $\xi_j := ( 2\tNF)^{\frac 1 2} (j\pi - \alpha)^{\frac 1 2 } $, i.e.\
\begin{equation}
\begin{aligned}
 B_j &= \{ \bxi \in \mR^m : |\bxi| \in (b_{j-}; \; b_{j+})\}, 
&&b_{j\pm} := (\xi_j^2 \pm 2\tNF \varepsilon )^{\frac 1 2}   \MTEXT{for} j \in \mN \\
 B_0 &= \{ \bxi \in \mR^m : |\bxi| < b_0 \},
&&b_0 := (2\tNF)^{\frac 1 2}  \max(  \varepsilon   - \alpha, 0)^{\frac 1 2}. 
\end{aligned}
\label{eq:IP2DefBj}
\end{equation}
Note that the constructed Fourier domain splitting is disjoint, i.e.\ $\mR^m =  ( \bigsqcup_{j=0}^\infty B_j ) \sqcup D_\varepsilon$.
From the definition of the domain $D_{\varepsilon}$, it immediately follows that
\begin{equation}
J_\varepsilon \geq \sin(\varepsilon) \norm{\hat \varphi|_{D_\varepsilon}}. \label{eq:IP2CTFBulkBound}
\end{equation} 
Hence, what remains to be done is to derive bounds for the sub-integrals $J_j$ around the zero-manifolds of $s_\alpha$ and to choose $\varepsilon$ to balance the contributions in \eqref{eq:IP2StabDecomp} . 

% We achieve this by approximating  $s_\alpha$ locally as a polynomial and estimating the arisig integrals via the following Lemma:

 \vspace{.5em}
\subsection{Estimate for the central CTF-minimum} \label{SS4.2}

We first consider the ball-shaped domain $B_0$, i.e.\ the low frequency part of 
the Fourier domain splitting.
Note that $B_0$ does not contain a zero of $s_\alpha$ if $\alpha > 0$ but still a (possibly small) local minimum of $s_\alpha^2$ at $\bxi = 0$.
Since $B_0 = \emptyset$ for $\alpha \geq \varepsilon$, we may restrict to the case $ \alpha < \varepsilon$.

By \eqref{eq:IP2DefBj}  and the assumption $\varepsilon\leq \pi/6$, we have $|\bxi|^2/(2\tNF)   + \alpha \in  [0;\pi/6)$ and hence $s_\alpha(\bxi)^2 \geq C_0 (|\bxi|^2/(2\tNF)   + \alpha)^2$ for all $\bxi \in B_0$ with $C_0:= \sin(\pi/6)^2/(\pi/6)^2$. This implies
\begin{equation}
\begin{aligned}
 J_0 &= \int_{B_0}  s_\alpha(\bxi)^2  |\hat \varphi ( \bxi) |^2  \; \D\bxi  \geq C_{0 } \int_{B_0}   \left( \frac{ |\bxi|^4 }{4\tNF^2}  +  \frac{   \alpha |\bxi|^2 }{\tNF} +  \alpha^2 \right)  |\hat \varphi ( \bxi) |^2  \; \D\bxi 
 \end{aligned}
 \label{eq:IP2Order2Zero1}
\end{equation}
From \eqref{eq:IP2Order2Zero1}, it can be seen that $J_0$ is bounded from below by $\alpha^2 \norm{\hat \varphi|_{B_0}}^2$ and thus indeed stable if $\alpha \neq 0$. However, since $\alpha$ is typically (almost) zero in hard X-ray imaging, we have to resort to the integral summands that involve powers of $|\bxi|^2$ in order to achieve robust estimates. These integrands have a zero at $\bxi = 0 \in B_0$, which is why smoothness of $|\hat \varphi |^2$ has to be exploited to obtain reasonable estimates. This is achieved by the following lemma:

\vspace{.5em}
\begin{lemma} \label{lem:IntegralDer2Bound}
Let $D:= \{ \bx \in \mR^m: |\bx| \leq a\}$ be a closed concentric ball of radius $a$ in $\mR^m$. Let $g,w \in \sC^2(  D)$  where $w$ is radially symmetric, i.e.\ $w(\bx) = w_0(|\bx|) $ for all $ \bx \in D$ and some function $w_0: [0;a] \to \mR$.
Then 
\begin{eqnarray*}
 \int_{D} \Delta w\cdot g \; \D\bx &=& \frac{mw_0'(a)}{a}  \int_{D} g \; \D\bx  
   + \int_{D} \left( \frac{w_0'(a)}{2a} (a^2 -  |\bx|^2) - (w_0(a)- w) \right)  \Delta g  \; \D\bx. %\label{eq:IP2Order2ZeroAux1}
\end{eqnarray*}
\end{lemma}
\vspace{.5em}
\begin{proof}
 By Green's second identity we have 
\begin{equation}
\begin{aligned}
\int_{D} \Delta w \cdot g \; \D\bx &= \int_{D} w \cdot \Delta g \; \D\bx + \int_{\partial D} \left( g \frac{\partial w}{\partial \bn} - w \frac{\partial g}{\partial \bn} \right) \; \D S(x)  \\
&= \int_{D} w \cdot \Delta g \; \D\bx + w_0'(a) \int_{\partial D} g \; \D S(x) - w_0(a) \int_{\partial D} \frac{\partial g}{\partial \bn} \; \D S(x). 
 \end{aligned}
  \label{eq:IntegralDer2Bound-1}
\end{equation}
Here, $\partial / ( \partial \bn )$ denotes the derivative along the unit normal vector pointing to the outside of $\partial D$.
 The boundary terms  can be eliminated via the relations
\begin{eqnarray*}
\int_{\partial D} \frac{\partial g}{\partial \bn} \; \D S(x) &=& \int_{D} \Delta g \; \D\bx \\
 \int_{\partial D} g  \; \D S(x) &=& \frac 1 {2a} \int_{\partial D}  g \frac{\partial (|\bx|^2) }{\partial \bn}  \; \D S(x) \\
 &=& \frac 1 {2a} \left( \int_{\partial D}  |\bx|^2 \frac{\partial g}{\partial \bn} \; \D S(x) + 2m\int_{D} g \; \D\bx - \int_{D} |\bx|^2 \Delta g \; \D\bx \right) \\
 &=& \frac 1 {2a} \left(  2m \int_{D} g \; \D\bx + \int_{D} (a^2 -  |\bx|^2) \Delta g \; \D\bx \right),
\end{eqnarray*}
which again follow from Green's second identity. Plugging this into \eqref{eq:IntegralDer2Bound-1} yields the claimed identity.   
\end{proof}
\vspace{.5em}

As the functions $\bxi \mapsto |\bxi|^2$ and $\bxi \mapsto |\bxi|^4$ are radially symmetric  and since $B_0$ is a concentric ball of radius $b_0$, we find that \cref{lem:IntegralDer2Bound} can be applied to the integrals in \eqref{eq:IP2Order2Zero1}. This yields the following estimate for the considered low-frequency subdomain:

\vspace{.5em}
\begin{lemma}%[Stability estimate for the central CTF-minimum]
\label{lem:IP2CTFCentralBound}
Let $\alpha < \varepsilon \leq \pi/6$ and  $B_0^{-} := \{\bxi \in B_0: \Delta | \hat \varphi |^2 (\bxi) \leq 0\}$. Then 
\begin{equation} \label{eq:IP2CTFCentralBound}
\begin{aligned}
 J_0 &\geq  C_0 \Big( \left( \tfrac{m}{m+4} (\varepsilon-\alpha)^2   
+ \tfrac{2m}{m+2}  \alpha (\varepsilon -\alpha)  +  \alpha^2  \right) 
\norm{\hat \varphi|_{B_0}}^2  \\
 &\;\;\;\;\;\;\;\;\; +  \tNF \left( \tfrac{2}{3(m+4)}\left( \varepsilon - \alpha \right) ^3  
+ \tfrac{1}{m+2}  \alpha \left( \varepsilon - \alpha \right)^2  \right)  \int_{B_0^-} \Delta | \hat \varphi|^2 \; \D\bxi \Big).
\end{aligned}
\end{equation}
\end{lemma}
\vspace{.5em}
\begin{proof}
For the functions $w_1(\bxi) :=  |\bxi|^6 / ( 6(m+4))$ and $w_2(\bxi) :=  |\bxi|^4 / ( 4(m+2))$, we have that $\Delta w_1 (\bxi) = |\bxi|^4$ and  $\Delta w_2 (\bxi) = |\bxi|^2$ for all $\bxi \in \mR^m$. Accordingly, the integral in \eqref{eq:IP2Order2Zero1} matches the setting of   \cref{lem:IntegralDer2Bound}. An application gives
 \begin{eqnarray*} 
  \lefteqn{\int_{B_0}   |\bxi|^4  |\hat \varphi (\bxi)|^2  \; \D\bxi =   \int_{B_0}   \Delta w _1(\bxi)   |\hat \varphi (\bxi)|^2  \; \D\bxi}  \\
\qquad&=& \frac{m b_0^4 }{m+4} \int_{B_0}     |\hat \varphi (\bxi)|^2  \; \D\bxi 
+  \frac{1 }{  m+4} \int_{B_0}  \left( \frac{b_0^4} 2 (b_0^2 - |\bxi|^2) - \frac{1} 6 (b_0^6 - |\bxi|^6) \right)   \Delta|\hat \varphi (\bxi)|^2  \; \D\bxi \\ 
&\geq&  \frac{m b_0^4 }{m+4} \int_{B_0}     |\hat \varphi (\bxi)|^2  \; \D\bxi  +  \frac{b_0^6 }{ 3( m+4)} \int_{B_0^{-}}   \Delta|\hat \varphi (\bxi)|^2.  \; \D\bxi   % \label{eq:IP2Order2Zero2-a}
\end{eqnarray*}
Here, the integral over $\Delta |\hat \varphi|^2 $ has been bounded from below by the integral within the sub-domain $B_0^{-} = \{\bx \in B_0 : \Delta |\hat \varphi|^2(\bx) \leq 0 \}$ multiplied by the 
maximum of the non-negative factor $(b_0^4/2) (b_0^2 - |\bxi|^2) - \frac{1} 6 (b_0^6 - |\bxi|^6)$. Analogously we obtain 
\begin{eqnarray*}
\lefteqn{\int_{B_0}   |\bxi|^2  |\hat \varphi (\bxi)|^2  \; \D\bxi = \int_{B_0}   \Delta w _2(\bxi)   |\hat \varphi (\bxi)|^2  \; \D\bxi }\\
&=& \frac{m b_0^2 }{m+2} \int_{B_0}     |\hat \varphi (\bxi)|^2  \; \D\bxi 
+  \frac{1 }{  m+2} \int_{B_0} \left( \frac{b_0^2} 2 (b_0^2 - |\bxi|^2) - \frac{1} 4 (b_0^2 - |\bxi|^2) \right)    \Delta|\hat \varphi (\bxi)|^2  \; \D\bxi \\
&\geq&  \frac{m b_0^2 }{m+2} \int_{B_0}     |\hat \varphi (\bxi)|^2  \; \D\bxi  +  \frac{b_0^4 }{ 4( m+2)} \int_{B_0^{-}}    \Delta|\hat \varphi (\bxi)|^2  \; \D\bxi   % \label{eq:IP2Order2Zero2-b}
\end{eqnarray*}
 Inserting these relations into \eqref{eq:IP2Order2Zero1} yields
\begin{eqnarray*}
J_0 &\geq& C_{0} \int_{B_0}   \left( \frac{ |\bxi|^4 }{4\tNF^2}  +  \frac{   \alpha |\bxi|^2 }{\tNF} +  \alpha^2 \right)  |\hat \varphi ( \bxi) |^2  \; \D\bxi  \\
&\geq& C_{0}   \left(\frac{m b_0^4 }{4\tNF^2(m+4)} +  \frac{m \alpha b_0^2}{\tNF( m+2)}  + \alpha^2 \right)  \int_{B_0}    |\hat \varphi ( \bxi) |^2  \; \D\bxi \\ 
&+& C_{0}\left(\frac{ b_0^6 }{12 \tNF^2 (m+4)}+ \frac{ b_0^4 \alpha }{4\tNF(m+2)} \right) \int_{B_0^{-}}    \Delta  |\hat \varphi ( \bxi) |^2 \; \D\bxi.  
\end{eqnarray*} 
Re-substituting $b_0 = (2\tNF)^{\frac 1 2 }  (\epsilon - \alpha )^{ \frac 1 2 }$ into this estimate gives the assertion.  % (recall that $\varepsilon > \alpha$) 
\end{proof}

\vspace{.5em}

\subsection{Estimate for the first order CTF-zeros} \label{SS4.3}

Before we can derive a global stability estimate from \eqref{eq:IP2CTFCentralBound}, we have to consider the remaining integrals $J_j$, $j \geq 1$ in the Fourier domain 
splitting \eqref{eq:IP2StabDecomp}.
The integration domains $B_j$ are annular shells around the first order zero-manifolds $Z_j = \{ \bxi \in \mR^m: |\bxi| = \xi_j\}$,  $j\in \mN$ of the CTF $s_\alpha$.
According to \eqref{eq:IP2DefBj}, we have $|\bxi|^2/(2\tNF) + \alpha - j\pi \in  (-\varepsilon;  \varepsilon)$ for all $\bxi \in B_j$ and thus
\begin{equation}
 s_\alpha(\bxi)^2 = \sin\left( \frac{|\bxi|^2}{2\tNF} \!+\! \alpha  \!+\! j \pi \right)^2 \geq \frac{\sin(\varepsilon)^2 }{\varepsilon^2} \left( \frac{|\bxi|^2}{2\tNF} \!+\! \alpha  \!+\! j \pi \right)^2 \stackrel{\varepsilon \leq \frac \pi 6}\geq   \frac{C_0}{ 4 \tNF^2}  \left( |\bxi|^2 - \xi_j^2 \right)^2.
\end{equation}
Using that $(|\bxi|^2 - \xi_j^2)^2   = (|\bxi| + \xi_j)^2 (|\bxi| - \xi_j)^2 \geq (b_{j-} + \xi_j)^2 (|\bxi| - \xi_j)^2 $, this yields
\begin{equation}
 J_j = \int_{B_j}  s_\alpha(\bxi)^2   |\hat \varphi ( \bxi) |^2\D\bxi 
\geq   \frac{ C_0  (b_{j-}\! +\! \xi_j)^2  }{4 \tNF^2} \tilde{J}_j,\quad 
\tilde{J}_j:= \int_{B_j} \left( |\bxi| - \xi_j \right)^2  |\hat \varphi ( \bxi) |^2   \; \D\bxi. \label{eq:IP2Order1Zero1}
\end{equation}
for $j \in \mN$.
Accordingly, we have to estimate polynomially weighted integrals, similar to the preceding section, yet within the shell-shaped domains $B_j$ instead of a concentric ball. We achieve this by introducing polar coordinates and applying \cref{lem:IntegralDer2Bound} to the resulting one-dimensional radial integrals:

\vspace{.5em}
\begin{lemma}[Stability estimate for the first order CTF-zeros] \label{lem:IP2CTF1stOrderBound}
Let $0 < \varepsilon \leq \pi/6$ and $j \in \mN$ . Then there exists a measurable subset $B_j^{-} \subset B_j$ and a constant $0 < C_1 \leq 1/4$, which depends only on $m$,
such that
\begin{equation}
 J_j \geq C_0 \left( C_1\varepsilon^2   \norm{\hat \varphi|_{B_j}}^2 + \frac{\varepsilon^4 \tNF}{ 32 ( j \pi - \alpha) } \int_{B_j^{-}} \Delta  |\hat \varphi|^2(\bxi) \; \D\bxi \right). \label{eq:IP2CTF1stOrderBound}
\end{equation}
\end{lemma}
\vspace{.5em}
\begin{proof}
 Let $j\in \mN$ arbitrary. In order to make \cref{lem:IntegralDer2Bound} applicable to the present setting, we rewrite the integral $\tilde J_j$ in \eqref{eq:IP2Order1Zero1} by introducing polar coordinates $\bxi = \xi \btheta$, $\xi \geq 0$, $\btheta \in \mS^{m-1}$ ($\mS^{m-1}$: unit sphere in $\mR^m$). With the notation in \eqref{eq:IP2DefBj}, this yields \vspace{-.6em}
 \begin{equation}
 \tilde J_j = \int_{b_{j-}}^{b_{j+}} (\xi - \xi _j)^2 \xi^{m-1} 
\varphi_{\textup{rad}}(\xi)\,\D\xi,\qquad 
\varphi_{\textup{rad}}(\xi):= \int_{\mS^{m-1}} |\hat \varphi ( \xi \btheta ) |^2  \; \D \theta . 
\end{equation}
Setting $a_j := b_{j+}- \xi_j$, $\eta := \xi - \xi_j$, 
$\varphi_{j}(\xi-\xi_j):=  \xi ^{(m-1)/2} \varphi_{\textup{rad}}(\xi)$,
%and using that $\xi_j - b_{j-} > a_j \geq 0$
we obtain the bound 
\begin{equation}
\begin{aligned} 
  \tilde J_j &= \int_{\xi_j - a_j}^{\xi_j + a_j} (\xi - \xi _j)^2 \xi^{\frac{m-1} 2} \varphi_{j}(\xi-\xi_j) \,\D\xi 
	+ \int^{\xi_j - a_j}_{b_{j-}} (\xi - \xi _j)^2\xi^{m-1} \varphi_{\textup{rad}}(\xi) \; \D\xi   \\
  &\geq (\xi_j-a_j)^{\frac{m-1} 2} \int_{-a_j}^{a_j} \eta^2 \varphi_j(\eta) \; \D\eta + a_j^2   \int^{\xi_j - a_j}_{b_{j-}} \xi^{m-1} \varphi_{\textup{rad}}(\xi) \; \D\xi. 
\end{aligned}  \label{eq:IP2Order1Zero1-1}
\end{equation}
The first integral on the right hand side of \eqref{eq:IP2Order1Zero1-1} matches the setting of \cref{lem:IntegralDer2Bound} in $m=1$ dimensions with weight function $w(\eta):= \eta^4/12$. This yields 
\begin{equation}
\begin{aligned} 
 \int_{-a_j}^{a_j}    \eta^2 \varphi_j(\eta) \; \D\eta &= \frac{a_j^2} 3  \int_{-a_j}^{a_j}  \varphi_j(\eta) \; \D\eta  
  +   \int_{-a_j}^{a_j} \left( \tfrac{a_j^2  (a_j^2 - \eta^2)} 6 - \tfrac{a_j^4 - \eta^4}{12} \right)\varphi_j''(\eta) \; \D\eta  \\
 &\geq  \frac{a_j^2} 3  \int_{-a_j}^{a_j}  \varphi_j(\eta) \; \D\eta + \frac{a_j^4}{12}  \int_{I^-_j}  \varphi_j''(\eta) \; \D\eta,
\end{aligned}  \label{eq:IP2Order1Zero1-2}
\end{equation}
where  we have set $I_j^{-}:=\{ \eta\in [-a_j; a_j] : \varphi_j''(\eta) \leq 0\}$ and used that
$0\leq (a_j^2/6)  (a_j^2 - \eta^2) - (a_j^4 - \eta^4) / 12 \leq a_j^4/12$.
Setting  $B_j^{-} := \{ \bxi \in B_j : |\bxi|   \in I_j^- + \xi_j \}$
and re-substituting $\varphi_{\textup{rad}}$ and $|\hat \varphi|^2$ for $\varphi_j$, we estimate the second integral in \eqref{eq:IP2Order1Zero1-2} by  
\begin{align}
\label{eq:IP2Order1Zero1-3} \lefteqn{\underbrace{(\xi_j-a_j)^{\frac{m-1} 2}}_{\leq \xi^{(m-1)/2} \text{ for } \xi   \in I_j^{-}+ \xi_j} \int_{I_j^- }  \underbrace{ \varphi_j''(\eta)}_{\leq 0} \; \D\eta \geq \int_{I_j^- + \xi_j}   \xi^{\frac{m-1} 2} \partial_{\xi}^2 \left(  \xi^{\frac{m-1} 2} \varphi_{\textup{rad}}(\xi) \right) \; \D\xi}   \\
&= \int_{I_j^- + \xi_j}   \xi^{m-1} \int_{\mS^{m-1}} \bigg( \partial_{\xi}^2  + \frac{m-1}{\xi} \partial_\xi + \frac{(m-1)(m-3)}{4\xi^2}\bigg) |\hat \varphi(\xi \btheta)|^2 \; \D \theta  \; \D\xi \nonumber \\
% = \int_{I_j^- + (k\pi - \varphi_0)^{\frac 1 2 }}  \xi^{\frac{m-1} 2} \int_{\mS^{m-1}} \frac{\partial^2}{\partial \xi^2} \left( \xi^{\frac{m-1} 2} \tilde f(\xi \btheta)  \right)\; \D \theta  \; \D\xi \\
&= \int_{B_j^{-}} \!\left( \Delta + \frac{(m-1)(m-3)}{4|\bxi|^2}  \right) |\hat \varphi|^2(\xi \btheta) \; \D\bxi   \nonumber \\
&\geq \int_{B_j^{-}} \!\!\Delta  |\hat \varphi|^2(\bxi) \; \D\bxi - \frac{\delta_{m2}}{4 b_{j-}^2 }  \int_{B_j^-}  \! |\hat \varphi|^2 (\bxi)  \; \D\bxi.  \nonumber
\end{align}
Here we have identified the radial part of the Laplacian in polar coordinates and implicitly added the angular part, which does not contribute to the $\mS^{m-1}$-integrals. The inequality in the last line follows from the fact that $(m-1)(m-3)/(4|\bxi|^2)$ is negative only for $m=2$, in which case it is bounded by the given term ($\delta_{ij}$: Kronecker-Delta). Re-substituting $|\hat \varphi|^2$ in the remaining integrals in \eqref{eq:IP2Order1Zero1-1} and \eqref{eq:IP2Order1Zero1-2} gives 
\begin{align}
  &\;\;\;\;\; (\xi_j-a_j)^{\frac{m-1} 2} \frac{a_j^2} 3  \int_{-a_j}^{a_j}  \varphi_j(\eta) \; \D\eta 
+ a_j^2   \int^{\xi_j - a_j}_{b_{j-}} \xi^{m-1} \varphi_{\textup{rad}}(\xi) \; \D\xi
\label{eq:IP2Order1Zero1-4} \\
 &\geq  \left( \tfrac{\xi_j-a_j}{\xi_j+a_j} \right) ^{\frac{m-1} 2}  \frac{a_j^2} 3  \int_{\xi_j-a_j}^{\xi_j+a_j}  \xi^{m-1}  \varphi_{\textup{rad}}(\xi) \; \D\xi  + a_j^2   \int^{\xi_j - a_j}_{b_{j-}} \xi^{m-1} \varphi_{\textup{rad}}(\xi) \; \D\bxi \nonumber \\
 &\geq  \left( \tfrac{\xi_j-a_j}{\xi_j+a_j} \right) ^{\frac{m-1} 2}  \frac{a_j^2} 3 \int_{b_{j-}}^{b_{j+}} \! \xi^{m-1}\!   \int_{\mS^{m-1}} \!  |\hat \varphi|^2(\xi\btheta) \; \D \theta \; \D\bxi =\left( \tfrac{\xi_j-a_j}{\xi_j+a_j} \right) ^{\frac{m-1} 2}  \frac{a_j^2} 3 \norm{\hat \varphi|_{B_j}}^2. \nonumber
\end{align}
Combining the estimates \eqref{eq:IP2Order1Zero1-1}, \eqref{eq:IP2Order1Zero1-2}, \eqref{eq:IP2Order1Zero1-3} and \eqref{eq:IP2Order1Zero1-4}, we thus obtain for all $j \in \mN$
\begin{eqnarray}
 \tilde J_j \geq   \underbrace{\frac{a_j^2} 3   \left(  \left( \frac{\xi_j-a_j}{\xi_j+a_j} \right) ^{\frac{m-1} 2}   -  \frac{a_j^2}{16 b_{j-}^2 } \right)}_{=:\lambda}   \norm{\hat \varphi|_{B_j}}^2 + \underbrace{\frac{a_j^4}{12}}_{=:\mu} \int_{B_j^{-}} \Delta  |\hat \varphi|^2(\bxi) \; \D\bxi  \label{eq:IP2Order1Zero1-5}
\end{eqnarray}

What remains to be done is to derive uniform bounds for the constants in \eqref{eq:IP2Order1Zero1-5} within the assumed parameter range $j \in \mN$,  $\alpha \leq \pi /2$ and $0 < \varepsilon \leq \pi/6$. 
First, we have $\xi_j^2 =  2\tNF  ( \pi j - \alpha) \geq \pi \tNF$ and hence by Taylor-expansion of the square-root
\begin{equation}
 a_j  = (\xi_j^2 + 2\tNF\varepsilon)^{\frac 1 2} -\xi_j  \in \frac{\varepsilon \tNF }{\xi_j} \cdot \left[1 - \frac{\varepsilon \tNF}{2 \xi_j^2 } ; 1 \right]  \subset \frac{\varepsilon \tNF}{  \xi_j} \cdot \left[ 1- \frac{\varepsilon}{2\pi} ; 1 \right]  \subset \frac{\varepsilon \tNF}{  \xi_j} \cdot \left[ \frac{11 }{12}; 1 \right] . \label{eq:IP2Order2ZeroAux3-1}
\end{equation}
% Moreover, we see that $\xi_j a_j / (\varepsilon \tNF)$ tends to one (uniformly in $j$) for $\varepsilon \to 0$.
Using \eqref{eq:IP2Order2ZeroAux3-1} along with $b_{j-}^2 = \xi^2_j - 2\tNF \varepsilon \stackrel{\varepsilon\leq \pi /6}\geq  \xi_j^2 - \pi \tNF / 3 \stackrel{\xi_j^2 \geq \pi \tNF}\geq  2\pi \tNF / 3$ 
furthermore gives
\begin{subequations} \label{eq:IP2Order2ZeroAux3-2} 
 \begin{align}
%  &&a_j  = (\xi_j^2 + 2\tNF\varepsilon)^{\frac 1 2} -\xi_j  \in \frac{\varepsilon}{\xi_j} \cdot \left[1 - \frac{\varepsilon \tNF}{ \xi_j^2 } ; 1 \right]  \subset \frac{\varepsilon}{2 \xi_j} \cdot \left[ \frac 5 6 ; 1 \right]\label{eq:IP2Order2ZeroAux3-1} \\
 \frac{a_j^2}{16 b_{j-}^2}  \stackrel{\text{\eqref{eq:IP2Order2ZeroAux3-1}}}\leq   \frac{\varepsilon ^2 \tNF^2}{16\xi_j^2 (\xi_j^2 - 2\tNF \varepsilon)} \leq \frac{3  \varepsilon^2 }{32\pi^2 } \leq \frac{1}{384}  \\
 \left( \frac{\xi_j-a_j}{\xi_j+a_j} \right)^{\frac 1 2 } \geq  1 - \frac { a_j }{\xi_j}   \stackrel{\text{\eqref{eq:IP2Order2ZeroAux3-1}}}\geq     1 - \frac { \varepsilon \tNF }{\xi_j^2}   \geq  1 - \frac { \varepsilon  }{\pi} \geq \frac 5 6
\end{align}
\end{subequations}     
Combining  \eqref{eq:IP2Order2ZeroAux3-1} and \eqref{eq:IP2Order2ZeroAux3-2} gives $ \lambda \geq (4/3) C_1 \varepsilon^2 \tNF^2  / \xi_j^2 $ with $C_1 :=  (121/576) ((5/6)^{m-1}  - \delta_{m2}/384)$ and $\mu \leq \varepsilon^4 \tNF^4 / ( 12 \xi_j^4)$. Substituting these bounds into \eqref{eq:IP2Order1Zero1-5}, we obtain
\begin{equation}
 \tilde J_j \geq \frac{4\tNF^2}{ 3 \xi_j^2}  \left(  C_1 \norm{\hat \varphi|_{B_j}}^2 + \frac{\varepsilon^4 \tNF^2}{ 16 \xi_j^2}  
\int_{B_j^{-}} \Delta  |\hat \varphi|^2(\bxi) \; \D\bxi \right). \label{eq:IP2Order1Zero1-6}
\end{equation}
Now the assertion follows by inserting \eqref{eq:IP2Order1Zero1-5} into \eqref{eq:IP2Order1Zero1}, using that  $\tNF / \xi_j^2 = 1/ (2(j \pi - \alpha))$ and $(\xi_j + b_{j-})^2 > 3\xi_j^2$ according to the estimate $b_{j-}^2 \geq 2\xi_j^2/3$.
\end{proof}
\vspace{.5em}
% 
% 
% We emphasize that the undesirable, yet unavoidable second summand in the estimates \eqref{eq:IP2CTF1stOrderBound} decays like $\cO( j^{-1})$ with the index $j$ of the considered zero-manifold $Z_j \subset B_j$. This means that zeros $\xi_j$ of the CTF $s_\alpha$ with large $j$ are less critical for the stability of \cref{ip2} under the assumed support constraint. This can be understood from \cref{figure5}: the gradient $|\nabla s_\alpha(\bxi)|$ at the zero-manifold $\bxi \in Z_j$ increases with $j$ so that - owing to the finitely sharp peaks of the object's Fourier transform $\hat \varphi$ - the obtained contrast $\norm{s_\alpha \cdot \hat \varphi|_{B_j}}$ is in g larger obtained for larger $j$.

\subsection{Global stability results} \label{SS4.4}

With the bounds for different subdomains of the Fourier space obtained in \cref{lem:IP2CTFCentralBound,lem:IP2CTF1stOrderBound}, we are now in a position to prove a global stability result for \cref{ip2}. For $\varepsilon > \alpha$, an application of the estimates \eqref{eq:IP2CTFBulkBound}, \eqref{eq:IP2CTFCentralBound} and \eqref{eq:IP2CTF1stOrderBound} to the sub-integrals in \eqref{eq:IP2StabDecomp} 
and \cref{lem:CpctSuppDer2Bound} with $\Omega = \closure{B(0,1/2)}$ yields
\begin{equation}\label{eq:IP2StabGlobal}
\begin{aligned} 
\norm{ s_\alpha \cdot \hat \varphi }^2  &\geq C_0 \left(  \zeta(\varepsilon, \alpha) \norm{ \varphi}^2 + \tNF \eta(\varepsilon, \alpha) \int_{B^-} \Delta | \hat \varphi|^2 \; \D\bxi \right) \\ 
&\geq  C_0 \left(  \zeta(\varepsilon, \alpha)  - \tfrac{\tNF} 2 \eta(\varepsilon, \alpha)  \right) \norm{\varphi}^2 
\end{aligned}
\end{equation}
Here, we have used that $\norm{\hat\varphi|_{D_\varepsilon}}^2 + \sum_{j=0}^\infty  \norm{\hat\varphi|_{B_j}}^2 = \norm{\hat \varphi}^2 = \norm{\varphi}$ by the disjoint Fourier 
 domain splitting in \secref{SS4.1} and defined $B^-:=  \bigcup_{j=0}^\infty B_j^-$ as well as
\begin{subequations}
\begin{equation}
 \begin{aligned} \zeta(\varepsilon, \alpha)  &:= 
\min \left\{ \tfrac{m}{m+4} (\varepsilon-\alpha)^2   + \tfrac{2m}{m+2}  \alpha (\varepsilon -\alpha)  +  \alpha^2 , C_1 \varepsilon^2  \right\}, \\
  \eta(\varepsilon, \alpha)  &:=  \max\left\{ 
	\tfrac{2}{3(m+4)} \left( \varepsilon - \alpha \right) ^3  + \tfrac{1}{m+2} \alpha \left( \varepsilon - \alpha \right)^2 , \tfrac{1}{ 32 (   \pi - \alpha) } \varepsilon^4  \right\} 
\end{aligned} 
\;\;\text{ for }\;\; \varepsilon > \alpha. \label{eq:IP2StabConstants-1}
 \end{equation}
Note that the bounding constant $\sin(\varepsilon)$ for the integral $J_{\varepsilon}$ does not appear in $\zeta(\varepsilon, \alpha)$ as $\sin(\varepsilon) \geq C_1 \varepsilon^2$ holds true within the entire parameter range $0\leq \varepsilon \leq \pi / 6$.

In the case $\varepsilon \leq \alpha$, the subdomains $B_0^- \subset B_0 = \emptyset$ are empty so that contributions from the sub-integral $J_0$ can be suppressed in the  estimate \eqref{eq:IP2StabGlobal}. Hence, we may set 
\begin{equation}
 \zeta(\varepsilon,\alpha)  :=  C_1 \varepsilon^2, \quad \eta(\varepsilon, \alpha) := \frac{\varepsilon^4 }{ 32 (   \pi - \alpha) }
 \;\;\text{ for }\;\; \varepsilon \leq \alpha. \label{eq:IP2StabConstants-2}
\end{equation}
\end{subequations}

Since $\eta(\varepsilon, \alpha)$ is of higher-order in $\varepsilon$  than $\zeta(\varepsilon, \alpha)$ according to \eqref{eq:IP2StabConstants-1} and \eqref{eq:IP2StabConstants-2}, it is possible to find an optimal value $\varepsilon = \varepsilon_{\textup{opt}}(\alpha, \tNF)$ such that the constant on the right-hand side of \eqref{eq:IP2StabGlobal} is maximal and in particular positive. This idea leads to the sought estimates for the stability constant $C_{\textup{IP2}}(\Omega,\tNF, \alpha)$ in \cref{thm:IP2StabRes}:

\vspace{.5em}
 \begin{proof}[Proof of \cref{thm:IP2StabRes}.]
 We first study the regime $\varepsilon \geq \alpha$, as is necessary in particular if $\alpha = 0$. Since 
 \begin{eqnarray*}
  \lefteqn{\tfrac{m}{m+4} (\varepsilon-\alpha)^2  + \tfrac{2m}{m+2}  \alpha (\varepsilon -\alpha)  +  \alpha^2 }\\ 
 & =& \tfrac{m}{m+4} \varepsilon^2 + (\tfrac{2m}{m+2}-\tfrac{2m}{m+4})\alpha(\varepsilon - \alpha) + (\tfrac{m}{m+4} - \tfrac{2m}{m+2} + 1) \alpha^2 
\geq \tfrac{m}{m+4} \varepsilon^2,
 \end{eqnarray*}
 \eqref{eq:IP2StabConstants-1} yields the bound %$\zeta(\varepsilon, \alpha) \geq \zeta_0 \varepsilon^2$
  \begin{equation} \label{eq:IP2StabRes-Proof-1}
  \zeta(\varepsilon, \alpha) \geq \zeta_0 \varepsilon^2, \quad \zeta_0:= \min\{ C_1, \tfrac{m}{m+4} \}.
  \end{equation}
% where $\zeta_0 :=  C_1 /5 $ in the case $m=1$ and 
% $\zeta_0 := C_1 / 3$ for $m>1$. 
To bound $\eta$ in \eqref{eq:IP2StabConstants-1} we observe that 
$\alpha \left( \varepsilon - \alpha \right)^2\leq 4 \varepsilon^3 / 27$ for 
$0 \leq \alpha < \varepsilon$ to obtain
 \[
  \tfrac{2}{3(m+4)} \left( \varepsilon - \alpha \right) ^3  + \tfrac{1}{m+2}\alpha \left( \varepsilon - \alpha \right)^2  \leq (\tfrac{2}{3(m+4)}+\tfrac{4}{27(m+2)}) \varepsilon^3
	\leq (\tfrac{2}{15}+\tfrac{4}{81})\varepsilon^3 < \tfrac{1}{5}\varepsilon^3.
 \]
 As  $\tfrac{\varepsilon  }{ 32(\pi - \alpha)}<\tfrac{1}{5}$ for 
$0\leq \alpha\leq\varepsilon\leq \tfrac{\pi}{6}$ this implies by \eqref{eq:IP2StabConstants-1}
 \begin{equation}
  \eta(\varepsilon, \alpha) \leq \tfrac{1}{5}\varepsilon^3.\label{eq:IP2StabRes-Proof-2}
 \end{equation}
 A comparison with \eqref{eq:IP2StabConstants-2} shows that the estimates \eqref{eq:IP2StabRes-Proof-1} and \eqref{eq:IP2StabRes-Proof-2} remain valid for $\varepsilon < \alpha$ and thus hold for all $0 \leq \varepsilon \leq \pi / 6$.

Now we apply the derived bounds on $\zeta(\varepsilon, \alpha)$ and $\eta(\varepsilon, \alpha)$ in \eqref{eq:IP2StabGlobal} to estimate the stability constant of 
\cref{ip2}:
\begin{align}\label{eq:IP2StabRes-Proof-3}
C_{\textup{IP2}} ( \Omega,\tNF, \alpha )^2 &= \inf_{\varphi \in L^2_\Omega, \, \norm{\varphi} = 1} \norm{S\varphi}^2 \stackrel{\eqref{eq:IP2Contrast}}=   
\inf_{\varphi \in L^2_\Omega, \, \norm{\varphi} = 1} 4  \norm{s_\alpha \cdot \hat \varphi }^2\nonumber\\
   &\geq
	\sup_{0\leq\varepsilon\leq \pi/6} 4 C_0 \varepsilon^2 \left( \zeta_0  - \tfrac{1}{10} \tNF \varepsilon  \right)    
  =   \begin{cases}  \frac{1600C_0\zeta_0^3}{27 \tNF^2}&\text{if }
	\tNF \geq \frac{40\zeta_0}{ \pi} \\
   \frac{\pi^2 C_0 }{9} ( \zeta_0 - \frac{\pi}{60} \tNF )   &\text{if }\tNF < \frac{40\zeta_0}{ \pi} \end{cases} \\
  &\geq  \min \left\{ c_1^2,  c_2^2 \tNF^{-2}  \right\}\nonumber
\end{align}
The supremum is attained at $\varepsilon = \varepsilon_{\textup{opt}} := \min\{ \pi /6, \, 20 \zeta_0 /(3\tNF) \} $, and we have defined $c_1:= (\pi^2C_0\zeta_0 / 27)^{1/2}$ and $c_2 :=( 1600C_0\zeta_0^3/ 27)^{1 / 2}$. 
%Note that we re-inserted the standard Fresnel number $\Fres = \Fres$ as used in the formulation of \cref{thm:IP2StabRes}.

It remains to show improved bounds for $\alpha>0$ and large $\tNF$. Here we choose 
$0 \leq \varepsilon \leq \alpha/3$ to satisfy $\varepsilon\leq \pi/6$ for all 
$\alpha\in (0,\pi/2]$. Moreover, we 
make use of the simplified expressions for $\zeta$ and $\eta$ in \eqref{eq:IP2StabConstants-2}, that apply in this parameter range. Inserting \eqref{eq:IP2StabConstants-2} into \eqref{eq:IP2StabGlobal} and using that $\pi - \alpha \geq \pi /2$ gives
 \begin{equation}
 \begin{aligned} 
   4\norm{ s_\alpha \cdot \hat \varphi }^2   &\geq 
	\sup_{0<\varepsilon\leq \alpha/3}
	4C_0 \varepsilon^2 \left( C_1  - \tfrac{1}{64(\pi - \alpha)} \tNF  \varepsilon^2  \right) \geq  
	\sup_{0<\varepsilon\leq \alpha/3}
	4C_0 \varepsilon^2 \left( C_1  -  \tfrac 1 {32\pi} \tNF \varepsilon^2  \right)  \\ 
   &=    \begin{cases} \frac{ 32 \pi C_0 C_1^2}{  \tNF}&\text{if }\tNF \geq \frac{144 \pi C_1}{ \alpha^2 } \\
   \frac{ 4C_0 \alpha^2  }{9} (C_1 - \frac{ \alpha^2 \tNF}{288 \pi}) &\text{if }\tNF < \frac{144 \pi C_1}{ \alpha^2 } \end{cases} \\
   &\geq \min \left\{  c_3^2  \alpha^2, c_4^2 \tNF^{-1}  \right\} 
\end{aligned} \label{eq:IP2StabRes-Proof-4}
\end{equation}
 for all $\varphi\in L^2_{\Omega}$ with  $\|\varphi\|=1$
 Here, the constants are chosen as $c_3:= (2C_0C_1/9) ^{1/2}$ and $c_4:=  (32\pi C_0)^{1/2} C_1$ and  the supremum is attained at $\varepsilon_{\textup{opt}} =   \min \{ \alpha/3 , \,  (16 \pi C_1 / \tNF )^{1/2}  \} $.
Combining \eqref{eq:IP2StabRes-Proof-3} and \eqref{eq:IP2StabRes-Proof-4} yields the claimed stability estimate \eqref{eq:IP2StabRes-2}. 
%Finally, by inserting the derived estimates \eqref{eq:IP2StabRes-Proof-3} and \eqref{eq:IP2StabRes-Proof-4} into \eqref{eq:IP2Contrast}, we obtain a characterization of the stability constant of \cref{ip2}:
%\begin{eqnarray*}
 %C_{\textup{IP2}} ( \Omega,\tNF, \alpha ) &=& \inf_{\varphi \in L^2_\Omega, \, \norm{\varphi} = 1} \norm{S\varphi} \stackrel{\eqref{eq:IP2Contrast}}=   \inf_{\cF^{-1}(\hat \varphi) \in L^2_\Omega, \, \norm{\hat \varphi} = 1} 2  \norm{s_\alpha \cdot \hat \varphi } \nonumber \\
 %&\geq& \max \left\{  \min \left\{ c_1, c_2 \Fres^{-1}  \right\}, \,  \min \left\{  c_3  \alpha , c_4  \Fres^{-1/2}  \right\} \right\}. \qquad \;\;\; \qedhere
%\end{eqnarray*}
\end{proof}
%\vspace{1em}

\vspace{.5em}

\subsection{Optimality} \label{SS4.5}

As opposed to the analysis of \cref{ip1}, the bounds on the stability constant $C_{\textup{IP2}}$ derived in this section can be expected to be highly non-optimal: the principal derivative-bound \eqref{eq:CpctSuppDer2Bound} is clearly not sharp and various additional estimates have been made in \cref{lem:IP2CTFCentralBound},  \cref{lem:IP2CTF1stOrderBound} and in the proof of \cref{thm:IP2StabRes} just to simplify the arising terms. Accordingly, the obtained numerical values for the constants $c_j$ in \eqref{eq:IP2StabRes-2} can be expected to be \emph{overly pessimistic}. Sharper characterizations may be achieved by numerically computing the minimal singular value of the forward operator $S$ by a similar approach as in \secref{SS3.5}.

Independently of such numerical refinements, it is of  interest whether the achieved analytical bounds $C_{\textup{IP2}} ( \Omega,\tNF, \alpha ) \gtrsim \tNF^{-\nu}$ 
%($\nu = 1$ for $\alpha = 0$, $\nu = 1/2$ otherwise)
are at least of \emph{optimal order} $\nu$ in the Fresnel number $\tNF$ for 
$\tNF \to \infty$. A positive answer is given by the following theorem:

\vspace{.5em}
\begin{theorem}[Order-optimality of the stability bounds for \cref{ip2}] \label{thm:IP2Optimality}
Within the setting of \cref{thm:IP2StabRes}, let $c'_1, c'_2, \nu > 0$ be constants such that
\begin{equation}
 C_{\textup{IP2}} ( \Omega,\tNF, 0 ) \geq \min\left\{ c_1', c_2' \tNF^{-\nu} \right\} \MTEXT{for all} \tNF >0. \label{eq:IP2Optimality}
\end{equation}
Then $\nu \geq 1$, i.e.\ the bound \eqref{eq:IP2StabRes-2} is of optimal asymptotic order in $\tNF$ for $\alpha = 0$.
\end{theorem}
\vspace{.5em}
\begin{proof}
 Assume $\nu < 1$. Let $\alpha = 0$ and $\varphi \in L^2_\Omega \cap \sC^2(\mR^m)$ with $\norm{\varphi} = 1$. Then we have by the definition of the forward operator in \eqref{eq:fwModel2}
 \begin{eqnarray}
  \norm{S_0\varphi}  &=& 2 \norm{s_\alpha \cdot \cF(\varphi)} 
	= 2 \left( \int_{\mR^m} \sin \left( \tfrac{|\bxi|^2 }{2\tNF} \right) ^2 |\cF(\varphi)(\bxi)|^2 \; \D\bxi \right)^{1/2}  \nonumber \\ 
  &\stackrel{\sin(x)^2 \leq x^2}\leq& \frac{1}{\tNF} \left( \int _{\mR^m} |\bxi|^2 |\cF(\varphi)(\bxi)|^2 \; \D\bxi  \right)^{1/2 } 
	= \tNF^{-1} \norm{\Delta \varphi }. \label{eq:IP2Optimality-proof1}
 \end{eqnarray}
 Note that $\norm{\Delta \varphi}\neq 0$ by the maximum principle since $\Omega$ is bounded. 
 Now let $\tNF > \max\{ (c'_1/c'_2)^{-1/\nu}, (2 \pi c_2' / \norm{\Delta \varphi })^{1/(1-\nu)} \}$. Then \eqref{eq:IP2Optimality} and \eqref{eq:IP2Optimality-proof1} imply
 \begin{eqnarray*}
  C_{\textup{IP2}} ( \Omega,\tNF, 0 ) &\stackrel{c_2'\tNF^{-\nu} < c_1'}=& c_2' \tNF^{-\nu} =  \frac{2 \pi c_2' \tNF^{1-\nu}}{\norm{\Delta \varphi }} \tNF^{-1} \norm{\Delta \varphi } 
  \stackrel{\eqref{eq:IP2Optimality-proof1}}\geq  \frac{2 \pi c_2' \tNF^{1-\nu}}{\norm{\Delta \varphi }} \norm{S_0\varphi}  > \norm{S_0\varphi}.
 \end{eqnarray*}
 Since $\norm{\varphi} = 1$, this contradicts the definition of the stability constant $C_{\textup{IP2}}$. Hence, $\nu \geq 1$ must hold by contradiction.
\end{proof}
\vspace{.5em}

Intuitively, the upper bound $-\nu \leq -1$ for the achievable order in \cref{thm:IP2Optimality} is related to the \emph{second} order zero of the CTF $s_\alpha$ at $\bxi = 0$ for $\alpha = 0$, which allows to bound $|s_{\alpha}|$ by a quadratic function. In the case $\alpha > 0$, this zero no longer exists which enables higher order behavior $C_{\textup{IP2}} \gtrsim \tNF^{-1/2}$ as shown in \cref{thm:IP2StabRes}. As the CTF still has \emph{first order} zero manifolds, we conjecture that this rate is likewise of optimal order. However, simple arguments as in the proof of \cref{thm:IP2Optimality} suffer from technical difficulties arising from the spherical geometry and the $\tNF$-dependence of the first order zero-manifolds. A treatment of the case $\alpha > 0$ is therefore omitted here.

\vspace{1em}

\section{Image reconstruction from two measurements} \label{S5}

% In this section, we consider the problem of image reconstruction from \emph{two} diffraction patterns, measured at different Fresnel numbers $\tNFa{1} \neq \tNFa{2}$. Our main concern is to show that stability of \cref{ip3} is closely related to properties of the forward operator $S$ of \cref{ip2}. Hence, the insights gained in the preceding section \secref{S4} directly carry over to the considered setting of image recovery from two holograms.

By axial translation of the object in \cref{figure1}(a), holograms $I_1,I_2$ may be recorded for different sample-detector distances and thus at different Fresnel numbers $\tNFa 1 \neq \tNFa 2$. It is often stated \cite{Cloetens1999,Jonas2004TwoMeasUniquePhaseRetr,Burvall2011TwoPlanes,Krenkel2014BCAandCTF} that the acquired additional data permits a more stable phase retrieval in this setting - in particular if phase shifts $\phi$ and attenuation $\mu$ are to be recovered as independent parameters. Within the weak object approximation (see \secref{SS2.2}), this  setting amounts to reconstructing $h = - \mu - \I \phi$ from measurements $(T_1 h , T_2 h)$, where $T_j$ denotes the linearized forward operator in \eqref{eq:fwModel1} to the Fresnel number $\tNF = \tNFa j$. Adopting the CTF-formulation in \eqref{eq:fwModelCTF}, the two-hologram setting is  modeled by the forward map
 \begin{align}\label{eq:fwModel3}
\begin{array}{ll}
  S^{(2)}: L^2 (\mR^m)^2 \!\! &\to L^2 (\mR^m)^2, \\ [1ex]
  \qquad \qquad \;\; \left(\begin{smallmatrix} \phi \\ \mu \end{smallmatrix}\right)  \!\! &\mapsto 
  2\cF^{-1} \big( \boldsymbol S \cdot  \cF 
	\left(\begin{smallmatrix} \phi \\ \mu \end{smallmatrix}\right)  \big),
	\end{array}
	\quad 
\boldsymbol S(\bxi):=\begin{pmatrix} \sin \left( \frac{|\bxi|^2 }{2\tNFa 1 } \right) & - \cos \left( \frac{|\bxi|^2 }{2\tNFa 1 } \right) \\
  \sin \left( \frac{|\bxi|^2 }{2\tNFa 2 } \right) & - \cos \left( \frac{|\bxi|^2 }{2\tNFa 2 }\right) \end{pmatrix}
 \end{align}
 Here, the matrix-vector-product is to be understood as a point-wise product of the vector- and matrix-valued function values and the (inverse) Fourier transforms $\cF, \cF^{-1}$ are meant component-wise.
 Furthermore, we denote by $\norm{f } := ( \norm{f_1}^2 + \norm{f_2}^2 )^{1/2}$ for $f = (f_1, f_2) \in L^2(\mR^m)^2$ the usual norm on $L^2(\mR^m)^2$. With the forward model defined in  \eqref{eq:fwModel3}, image recovery from two diffraction patterns may then be stated as:
 
 \vspace{.5em}
 \begin{MyIP}[Phase contrast imaging from two measurements]\label{ip3}
  For given support $\Omega \subset \mR^m$, recover the images  $\phi, \mu \in L^2_\Omega$ from two noisy holograms
  \begin{equation*} 
   \begin{pmatrix} 
    I^{\bepsilon}_1 -1 \\ I^{\bepsilon}_2-1
   \end{pmatrix}
  =   S^{(2)} \begin{pmatrix} 
    \phi \\ \mu
   \end{pmatrix} + \bepsilon \MTEXT{with}  \norm{\bepsilon}  \leq \epsilon.
   \end{equation*}
 \end{MyIP}
 \vspace{.5em}

As the same amount of object information is to be recovered from a larger data set, it is evident that \cref{ip3} cannot be more ill-posed or ill-conditioned than \cref{ip1}. In order to gain deeper insight, it is illustrative to write down an explicit inverse of the forward operator $S^{(2)}$ defined in \eqref{eq:fwModel3} by point-wise inversion of the contrast-transfer-matrix $\boldsymbol S$ as done in \cite{Gureyev2004CTF}: setting $\tNFa{-}:=  (\tNFa 1^{-1} - \tNFa 2^{-1} )^{-1}$ this yields for any $(\phi, \mu) \in L^2 (\mR^m )$
\begin{align}
 \begin{pmatrix} \phi \\ \mu \end{pmatrix}  &= \cF^{-1} \Bigg( \frac{1}{2 \sin\big(\frac{|\bxi|^2 }{2\tNFa{-}}\big)} \begin{pmatrix}  \cos \left( \frac{|\bxi|^2 }{2\tNFa 2 }\right) &  - \cos \left( \frac{|\bxi|^2 }{2\tNFa 1 } \right) \\
  \sin \left( \frac{|\bxi|^2 }{2\tNFa 2 } \right) &  - \sin \left( \frac{|\bxi|^2 }{2\tNFa 1 } \right)  \end{pmatrix} \cF S^{(2)} \begin{pmatrix} \phi \\ \mu \end{pmatrix} \Bigg). \label{eq:IP3-2}
\end{align}
All of the operations on the right-hand side of in \eqref{eq:IP3-2} are bounded except for the point-wise division by the factor $2 \sin\big(|\bxi|^2 / (2\tNFa{-}) \big)$. According to \eqref{eq:fwModel2}, we find that the latter operation exactly corresponds to solving \cref{ip2} for a pure phase object (i.e.\ $\alpha = 0$) at an effective Fresnel number of $\tNF = \tNFa{-}$. In particular, this implies that the solution of \cref{ip3} is unique but ill-posed without a support constraint, i.e.\ for general $\phi, \mu \in L^2(\mR^m)$.
By a more rigorous analysis of the analogy to \cref{ip2}, we obtain stability estimates:

 \vspace{.5em}
\begin{theorem}[Stability estimate for \cref{ip3}] \label{thm:IP3StabRes}
 Let $S_0$ denote the forward operator of \cref{ip2} for $\alpha = 0$ and $\tNF = \tNFa{-} :=  (\tNFa 1^{-1} - \tNFa 2^{-1} )^{-1}$. Then 
 \begin{align}
  \big \| S^{(2)} \left(\begin{smallmatrix} \phi \\ \mu \end{smallmatrix}\right) \big\| 
	\geq 2^{-\frac 1 2} \norm{ S_0 \left( \phi + \I\mu \right)}  \MTEXT{for all real-valued} \phi, \mu \in L^2 (\mR^m)  \label{eq:IP3StabRes-1}
 \end{align}
 In particular, for any support-domain $\Omega \subset \mR^m$, we have the relative stability estimate
  \begin{align}
 \big \| S^{(2)}\left( \begin{smallmatrix} \phi \\ \mu \end{smallmatrix} \right) \big \| \geq  2^{-\frac 1 2} C_{\textup{IP2}}(\Omega, \tNFa -, 0)  \norm{\left(  \begin{smallmatrix} \phi \\ \mu \end{smallmatrix} \right) } \MTEXT{for all real-valued} \phi, \mu \in L^2_\Omega.  \label{eq:IP3StabRes-2}
 \end{align} 
%  In particular, for the support-domain $\Omega := \closure{B(0,1/2)}$, we have the stability estimate
%   \begin{align}
%   \norm{ S^{(2)}\left( \begin{smallmatrix} \phi \\ \mu \end{smallmatrix} \right)} \geq  2^{-\frac 1 2} \min \{ c_1, c_2 \tNFa{-}^{-1}\}  \norm{\left(  \begin{smallmatrix} \phi \\ \mu \end{smallmatrix} \right) } \MTEXT{for real-valued} \phi, \mu \in L^2_\Omega,  \label{eq:IP3StabRes-2}
%  \end{align}
%  where $c_1, c_2 > 0$ denote the constants from the stability result in \cref{thm:IP2StabRes}.
\end{theorem}

 \vspace{.5em}

 \begin{proof}
 Let $\phi, \mu \in L^2 (\mR^m)$ be real-valued. Setting $\bhhat := (\cF(\phi), \cF(\mu))$ and exploiting that $\cF$ is unitary, we obtain by \eqref{eq:fwModel3}
 \begin{equation*}
  \big \| S^{(2)} \left(\begin{smallmatrix} \phi \\ \mu \end{smallmatrix}\right) \big \|^2 
	= 4 \big \| \boldsymbol S \bhhat \big \|^2 = 4 \int_{\mR^m} | \boldsymbol S(\bxi)  \bhhat(\bxi) |^2 \; \D\bxi. \label{eq:IP3-3}
 \end{equation*}
 The integrand  $| \boldsymbol S(\bxi)  \bhhat(\bxi) |^2$ is bounded from below by $\sigma_0^2(\bxi)|\bhhat(\bxi) |^2$ with the smallest singular value $\sigma_0(\bxi)$ of the 2-by-2 matrix $\boldsymbol S(\bxi)$. Direct computation shows
 \begin{equation*}
  \sigma_0^2(\bxi) = 1 - \left| \cos\left(  \tfrac{|\bxi|^2 }{2\tNFa{-}}  \right)  \right| \geq \frac 1 2 \sin\left(  \tfrac{|\bxi|^2 }{2\tNFa{-}}  \right)^2 \MTEXT{for all} \bxi \in \mR^m. \label{eq:IP3-4}
 \end{equation*}
 By inserting this result into the previous equation and comparing to \eqref{eq:fwModel2}, we obtain
  \begin{align*}
  \big \| S^{(2)} \left(\begin{smallmatrix} \phi \\ \mu \end{smallmatrix}\right) \big \|^2 
	&\geq 2 \int_{\mR^m} \sin\left(  \tfrac{|\bxi|^2 }{2\tNFa{-}}  \right)^2 |  \bhhat(\bxi) |^2 \; \D\bxi \nonumber \\
  &= 2 \int_{\mR^m} \sin\left(  \tfrac{|\bxi|^2 }{2\tNFa{-}}  \right)^2 | \cF(\phi)(\bxi) + \I \cF(\mu)(\bxi) |^2 \; \D\bxi = \frac 1 2 \norm{S_0(\phi + \I \mu)}^2. %\label{eq:IP3-5}
 \end{align*}
 This proves \eqref{eq:IP3StabRes-1}. The second inequality \eqref{eq:IP3StabRes-2} follows by bounding the right-hand side of \eqref{eq:IP3StabRes-1} with the stability constant of \cref{ip2}. %, as characterized in  \cref{thm:IP2StabRes} for ball-shaped support-domains.
\end{proof}
 \vspace{.5em}

We emphasize that \cref{thm:IP3StabRes} does not only relate the worst case stability of \cref{ip2,ip3} but \eqref{eq:IP3StabRes-1} identifies the contrast attained by individual modes under the forward operators $S_0$ and $S^{(2)}$. Specifically, whenever $\norm{S_0(\phi + \I \mu)}$ is large, the mode $(\phi, \mu)^{\textup T}$ also attains high contrast under the forward map $S^{(2)}$ of \cref{ip3}. It should furthermore be noted that the derived stability estimates are in terms of the \emph{difference Fresnel number} $\tNFa{-}^{-1} = \tNFa 1^{-1} - \tNFa 2^{-1} $. Hence, if the two holograms are recorded at similar Fresnel numbers $\tNFa 1 \approx \tNFa 2$, i.e.\ for only slightly varied setup parameters, $\tNFa{-}$ is very large so that the stability bounds in \cref{thm:IP3StabRes} will hardly be better than for the \emph{single} measurement case. This is quite intuitive as the second measurement  provides only little additional information in such a setting.

\vspace{1em}

\section{Discussion and conclusions} \label{S6}
In this paper we have studied the stability of phase retrieval in propagation-based phase contrast imaging within the linear contrast-transfer-function model (CTF) \cite{Guigay1977CTF,Turner2004FormulaWeakAbsSlowlyVarPhase}, valid for weakly interacting objects. While the image reconstruction problem is generally ill-posed and even severely non-unique if both phase-shifts $\phi$ and attenuation $\mu$ are to be recovered from a single (phaseless) intensity measurement (\cref{ip1}), we could show \emph{well-posedness} under a support constraint: if the image $h = -\mu - \I \phi$ is known to be supported in a bounded domain $\Omega \subset \mR^m$, then the attained contrast $\norm{Th} \geq C_{\textup{IP1}}(\Omega, \tNF) \norm{h}$ is bounded from below by a positive stability constant $C_{\textup{IP1}}(\Omega, \tNF) > 0$ according to \cref{thm:IP1StabRes}. 
%As established in \secref{S3}, this startling result is due to the fact that the inversion of $T$ may be related to a reconstruction problem from incomplete Fourier data, for which well-posedness results are readily available.

Numerical simulations (see \secref{SS3.5}) suggest that the proven lower bound  $C_{\textup{IP1}}(\Omega,  \tNF) \gtrsim \tNF^{1/4} \exp(-\tNF / 8)$ is quite sharp. Accordingly, $C_{\textup{IP1}}$ unfortunately decays nearly exponentially with the Fresnel number $\tNF$ of the imaging setup (computed w.r.t.\ the support diameter!), which translates into partly ridiculously small values within the typical range $10 \leq \Fres \leq 10^3$ in X-ray phase contrast experiments. For instance, the estimate \eqref{eq:IP1StabRes-2} gives $C_{\textup{IP1}}(\Omega,2\pi\cdot 100) \leq 10^{-33}$ at a Fresnel number $\Fres = 100$. In other words: although well-posed, the inversion of the forward operator $T$ is typically \emph{heavily ill-conditioned}, admitting amplification of data errors by bounded, yet \emph{huge} factors $C_{\textup{IP1}}^{-1}$ in the image reconstruction. It is thus not surprising that the independent recovery of both $\phi$ and $\mu$ from a single hologram is typically severely corrupted by artifacts in practice, see e.g.\ \cite{Ruhlandt2014,MaretzkeMasters}, and  commonly considered as not feasible \cite{Jonas2004TwoMeasUniquePhaseRetr,Nugent2007TwoPlanesPhaseVortex,Burvall2011TwoPlanes}.

On the other hand, such a reconstruction has been successfully demonstrated in \cite{MaretzkeEtAl2016OptExpr} up to slight low-frequency artifacts. In view of the present work, the key ingredient to this demonstration can be identified as the comparably small Fresnel number $\Fres \approx 14$ of the support in the considered setup. By the characterization of the least stable modes in \cref{thm:LeastStableModes}, this corresponds to a minimum contrast
\begin{equation}
 \norm{T \phi_{(0,0)} } \geq 10^{-4} \norm{\phi_{(0,0)} }  \MTEXT{and}  \norm{T \phi_{(0,1)} } \geq 7\cdot 10^{-4} \norm{\phi_{(0,1)} }
\end{equation}
for the least and second least stable image modes $\phi_{(0,0)} $ and $\phi_{(0,1)} $. The reconstruction of these low-frequency image components is thus feasible in principle, yet numerically cumbersome, which explains the residual artifacts in the reconstruction. 
Higher-order image modes attain contrast in the order of the data noise level or above and may thus be recovered with reasonable accuracy. Our stability estimates thus predict that the joint recovery of phase shifts $\phi$ and absorption $\mu$ is feasible if and only if $\Fres \lesssim 10$, i.e.\ in the \emph{deeply holographic regime} - in excellent agreement with numerical reconstructions.

%In order to establish stability also for larger Fresnel numbers $\tNF$ we studied two variants of \cref{ip1}: image recovery for homogeneous objects (including in particular non-absorbing objects), characterized by $\mu \propto \phi$, as well as reconstruction from \emph{two} holograms,  measured at different Fresnel numbers. For the corresponding \cref{ip2,ip3}, the obtained well-posedness results (see \cref{thm:IP2StabRes} and \cref{thm:IP3StabRes}) show much more favorable decay rates $\cO(\tNF^{-1})$ of the corresponding stability constants. The mathematical intuition behind these results is that the well-known problematic zeros of the contrast transfer functions are effectively regularized by the smoothness in Fourier space resulting from the assumed real-space support constraints, see \secref{SS4.1}. 
According to the improved stability estimate for the settings of \cref{ip2,ip3} 
(see \cref{thm:IP2StabRes,thm:IP3StabRes})
image recovery can be performed without additional regularization 
provided a sufficiently strong support constraint. This is consistent with 
the observed high performance of iterative projection algorithms in phase retrieval
(see \cite{BCL:02,Marchesini2007EvalPhaseRetrAlg,fienup:13}%Shechtman2015_PhaseRetrReview} 
and references therein, mostly for the 
far-field case and e.g.\ \cite{Giewekemeyer2011CellSketch,Bartels2012} for
near-field phase contrast), which make use of support constraints.
%\todo{die CTF-inversion methods werden hier jetzt nicht mehr erwähnt, da sie m.E.\ wegen der Erwähnung des far-field case nicht mehr passen.}
%, compared to direct CTF inversion methods, which do not incorporate support knowledge.

%In general, our results strongly highlight the considerable value of support constraints in propagation-based phase contrast imaging - establishing well-posedness for otherwise ill-posed and non-unique phase retrieval problems. However, 
Let us discuss some possible extensions of our results. First of all we have assumed that 
diffraction patterns can be measured in a whole detector plane, whereas real-world X-ray detectors necessarily cover only a finite area. If intensities are measured only in a bounded 
subdomain, the singular values of the forward operator $T$ will eventually decay 
exponentially since $T$ can be written as an integral operator with an analytic 
(though highly oscillitory) kernel. 
On the other hand, Fresnel propagation is mathematically equivalent to time-evolution within the free Schr\"odinger equation (if the $z$-coordinate in \Fref{figure1}(b) is identified with time). Localization properties of the latter model suggest that the finiteness
of the detection domain has little influence if it is chosen sufficiently large - in agreement with experimental and numerical experience. 
A better understanding and analysis of the impact of a finite field of 
view will be an interesting goal for future research.  
Further possible extensions include a treatment of the nonlinear imaging model \eqref{eq:fwModelNL}, 
non-plane wave illumination and partial coherence.  
%So far, however, the scope of our results is limited to the (frequently used) weak object linearization. Stability analysis of the general, nonlinear imaging model \eqref{eq:fwModelNL} will be subject of future work. Likewise, effects of various real-world deviations from the assumed idealized imaging model may be considered such as effects of finite field of views, non-plane wave illumination or limited coherence.
Moreover, in \emph{region-of-interest} imaging of extended objects, support constraints do not hold true so that stability needs to be established by other means.

Finally, we emphasize that - although seemingly specific to the considered imaging setup - our analysis treats a fairly general physical problem: the reconstruction of a (compactly supported) perturbation to a plane wave from intensities of the propagated field 
under the paraxial Helmholtz equation or - equivalently - within time-dependent Schr\"odinger dynamics. The results may thus be relevant for several related wave-optical and quantum-mechanical inverse problems.

\vspace{1em}

\section*{Acknowledgements}
We thank Tim Salditt, Aike Ruhlandt and Martin Kren-kel from the Institute for X-ray Physics at the University of G\"ottingen 
for enlightening discussions on phase contrast setups, a priori constraints and stability of image reconstruction,  as well as for providing the experimental data used 
in Figures \ref{figure1} and \ref{figure2} for visualization of the theory in this work.
%Financial support by the German Research Foundation DFG through Project C02 
%of the Collaborative Research Center 755 Nanoscale
%Photonic Imaging is likewise gratefully acknowledged.

\small
\bibliographystyle{siam}
\bibliography{./literature_NearfieldStability}

\end{document}